# REMOVABILITY OF THE FUNDAMENTAL SINGULARITY FOR THE HEAT EQUATION AND ITS CONSEQUENCES. I. KOLMOGOROV-PETROVSKY-TYPE TEST

UGUR G. ABDULLA

ABSTRACT. We prove the necessary and sufficient condition for the removability of the fundamental singularity, and equivalently for the unique solvability of the singular Dirichlet problem for the heat equation. In the measure-theoretical context the criterion determines whether the $h$-parabolic measure of the singularity point is null or positive. From the probabilistic point of view, the criterion presents an asymptotic law for conditional Brownian motion.

## 1. Prelude

Consider the fundamental solution of the heat equation:

$$(1.1) \qquad F(x,t) = \begin{cases} (4\pi t)^{-\frac{N}{2}} e^{-\frac{|x|^2}{4t}}, & x \in \mathbb{R}^N, t > 0, \\ 0, & x \in \mathbb{R}^N, x \neq y, t = 0. \end{cases}$$

It is a distributional solution of the Cauchy Problem

$$\mathcal{H}F := F_t - \Delta F = 0 \quad \text{in } \mathbb{R}^{N+1}_+, \ F = \delta \quad \text{on } \mathbb{R}^N \times \{t = 0\}$$

where $\delta$ is a unit measure with support at $x = 0$. For any fixed point $\gamma \in \mathbb{R}^N$, let

$$h(x,t) := F(x - \gamma, t)$$

be a fundamental solution with a pole at $\mathcal{O} := (\gamma, 0)$. Singularity of $h$ at $\mathcal{O}$ represents the natural phenomenon of the space-time distribution of the unit energy initially blown up at a single point. The fundamental singularity is non-removable for the heat equation in $\mathbb{R}^{N+1}_+$. In particular, the Cauchy Problem for the heat equation in $\mathbb{R}^{N+1}_+$ has infinitely many solutions in class $O(h)$.

The goal of this paper is to reveal the criterion for the removability of the fundamental singularity for open subsets of $\mathbb{R}^{N+1}_+$. Let $\Omega \subset \mathbb{R}^{N+1}_+$ be an arbitrary open set and $\partial\Omega \cap \{t = 0\} = \{\mathcal{O}\}$. Let $g : \partial\Omega \to \mathbb{R}$ be a boundary function such that $g/h$ is a bounded Borel measurable. Consider a *singular* **parabolic Dirichlet problem**(PDP):

$$(1.2) \qquad \mathcal{H}u = 0 \quad \text{in } \Omega, \ u = g \quad \text{on } \partial\Omega \setminus \{\mathcal{O}\}; \ u = O(h) \quad \text{at } \mathcal{O},$$

Without prescribing the behavior of $u/h$ at $\mathcal{O}$, there exists one and only one or infinitely many solutions of PDP (see Section 3.1, formula (2.8)). The main goal of this paper is to find a necessary and sufficient condition for open sets $\Omega$ for the uniqueness of the solution to the PDP without prescribing $u/h$ at $\mathcal{O}$. The problem of removability vs. non-removability of the fundamental singularity is equivalent to the question of the uniqueness of the solution to PDP (1.2) without prescribing the behavior of $u/h$ at $\mathcal{O}$. The following procedure provides a key problem in testing the removability of fundamental singularity. Let

$$\Omega_n := \Omega \cap \{t > n^{-1}\}, \ n = 1, 2, ...$$

and $u_n$ be a unique solution of the parabolic Dirichlet problem

(1.3) $\quad \mathcal{H}u = 0, \quad \text{in } \Omega_n; \ u|_{\partial\Omega_n \cap \{t > n^{-1}\}} = 0; \ u|_{\partial\Omega_n \cap \{t = n^{-1}\}} = h(x, n^{-1}).$

From the maximum principle it follows that

(1.4) $\quad 0 \leq u_{n+1}(x,t) \leq u_n(x,t) \leq h(x,t), \quad \text{on } \Omega_n.$

Therefore, there exists a limit function

(1.5) $\quad u_*(x,t) = \lim_{n \to +\infty} u_n(x,t), \ (x,t) \in \Omega,$

which satisfies (1.2), and

$$0 \leq u_*(x,t) \leq h(x,t), \ (x,t) \in \Omega.$$

The following is the key problem.

**Problem** $\mathcal{A}_\gamma$: *Is $u_* \equiv 0$ in $\Omega$ or $u_*(x,t) \not\equiv 0$ in $\Omega$? Equivalently, is fundamental singularity at $\mathcal{O}$ removable or nonremovable for $\Omega$?*

Next, we formulate the equivalent problem in $\mathbb{R}_-^{N+1}$. In that context, we are going to consider one-point Alexandrov compactification: $\mathbb{R}_-^{N+1} \to \mathbb{R}_-^{N+1} \cup \{\infty\}$. For any fixed finite $\gamma \in \mathbb{R}^N$, consider a function

(1.6) $\quad \tilde{h}(x,t) = e^{\langle x, \gamma \rangle + |\gamma|^2 t}.$

It is a positive solution of the heat equation in $\mathbb{R}_-^{N+1}$. If $\gamma = 0$, then $\tilde{h} \equiv 1$, while in the case when $\gamma \neq 0$, it is an unbounded solution with singularity at $\infty$. The singularity is not removable for the heat equation in $\mathbb{R}_-^{N+1}$. We aim to reveal the criterion for the removability of the fundamental singularity for the open subsets of $\mathbb{R}_-^{N+1}$.

Let $\tilde{\Omega} \subset \mathbb{R}_-^{N+1}$ be an arbitrary open set, and $g : \partial\tilde{\Omega} \to \mathbb{R}$ be a boundary function, such that $g/\tilde{h}$ is a bounded Borel measurable. Consider a *singular* **parabolic Dirichlet problem**(PDP):

(1.7) $\quad \mathcal{H}u = 0 \quad \text{in } \tilde{\Omega}, \ u = g \quad \text{on } \partial\tilde{\Omega}; \ u = O(\tilde{h}) \quad \text{at } \infty,$

Without prescribing the behavior of $u/\tilde{h}$ at $\infty$, there exists one and only one or infinitely many solutions of (1.7) (see Section 3.1, formula (2.15)). The alternation is equivalent to the question of removability vs. non-removability of the fundamental singularity at $\infty$ for $\tilde{\Omega}$. In particular, in the case $\gamma = 0$ ($\tilde{h} \equiv 1$), we are addressing the uniqueness of a bounded solution of (1.7), the problem solved in [1]. Similar to its counterpart (1.2), the key problem to test the removability of the singularity at $\infty$ is formulated as follows:

Let

$$\tilde{\Omega}_n := \tilde{\Omega} \cap \{t > -n\}, \ n = 1, 2, ...$$

and $\tilde{u}_n$ be a unique solution of the parabolic Dirichlet problem

(1.8) $\quad \mathcal{H}u = 0, \quad \text{in } \tilde{\Omega}_n; \ u|_{\partial\tilde{\Omega}_n \cap \{t > -n\}} = 0; \ u|_{\partial\tilde{\Omega}_n \cap \{t = -n\}} = \tilde{h}(x, -n).$

From the maximum principle it follows that

(1.9) $$0 \leq \tilde{u}_{n+1}(x,t) \leq \tilde{u}_n(x,t) \leq \tilde{h}(x,t), \quad \text{on } \tilde{\Omega}_n.$$

Therefore, there exists a limit function

(1.10) $$\tilde{u}_*(x,t) = \lim_{n \to +\infty} \tilde{u}_n(x,t), \ (x,t) \in \tilde{\Omega},$$

which satisfies (1.7), and

$$0 \leq \tilde{u}_*(x,t) \leq \tilde{h}(x,t), \ (x,t) \in \tilde{\Omega}.$$

The following is the key problem.

**Problem $\tilde{\mathcal{A}}_\gamma$**: Is $\tilde{u}_* \equiv 0$ in $\tilde{\Omega}$ or $\tilde{u}_*(x,t) \not\equiv 0$ in $\tilde{\Omega}$? Equivalently, is a fundamental singularity at $\infty$ removable or nonremovable for $\tilde{\Omega}$?

*Remark* 1.1. Problem $\tilde{\mathcal{A}}_\gamma$ can be formulated in $\mathbb{R}^{N+1}$ without any change. Indeed, the parabolic Dirichlet problem for the heat equation is uniquely solvable in any open subset $\mathbb{R}_+^{N+1}$ in a class $O(\tilde{h})$. Therefore, given arbitrary open set $\Omega \subset \mathbb{R}^{N+1}$, the solution of the parabolic Dirichlet problem in $\Omega$ can be constructed as a unique continuation of the solution to the parabolic Dirichlet problem in $\Omega_- = \Omega \cap \mathbb{R}_-^{N+1}$. Moreover, the latter is independent of the boundary values assigned on $\partial \Omega_- \cap \{t = 0\}$, since it is a parabolic measure null set for $\Omega_-$. This implies that the Problem $\tilde{\mathcal{A}}_\gamma$ is equivalent for $\Omega$ and $\Omega_-$. Otherwise speaking, the fundamental singularity at $\infty$ is removable for $\Omega \subset \mathbb{R}^{N+1}$ if and only if it is removable for $\Omega_-$.

The only problem in the family of formulated problems that is solved is the Problem $\tilde{\mathcal{A}}_\gamma$ when $\gamma = 0$ (or Problem $\tilde{\mathcal{A}}_0$). The Problem $\tilde{\mathcal{A}}_0$ was formulated by Kolmogorov in 1928 in the seminar on the probability theory at the Moscow State University in the particular case with $\Omega = \{|x| < f(t), -\infty < t < 0\} \subset \mathbb{R}^2$, with $f \in C(-\infty, 0]$ such that $f(-\infty) = +\infty, f \uparrow +\infty, (-t)^{-\frac{1}{2}} f \uparrow +\infty$ as $t \downarrow -\infty$. Kolmogorov's motivation for posing this problem was a connection to the probabilistic problem of finding asymptotic behavior at infinity of the standard Brownian path. Let $\{\xi(t) : t \geq 0, P_\bullet\}$ be a standard 1-dimensional Brownian motion and $P_\bullet(\mathbf{B})$ is the probability of the event $\mathbf{B}$ as a function of the starting point $\xi(0)$. Blumenthal's 01 law implies that $P_0\{\xi(t) < f(-t), t \uparrow +\infty\} = 0$ or 1; $f(-t)$ is said to belong to the upper class if this probability is 1 and to the lower class otherwise. Remarkably, Kolmogorov Problem's solution $u_*$ is $= 0$ or $> 0$ according to as $f(-t)$ is in lower or upper class accordingly. Kolmogorov Problem in a one-dimensional setting was solved by Petrovsky in 1935, and the celebrated result is called the Kolmogorov-Petrovski test in the probabilistic literature [27] (see also [4]).

The full solution of the Kolmogorov Problem for arbitrary open sets $\Omega$ (or Problem $\tilde{\mathcal{A}}_0$) is presented in [1]. A new concept of regularity or irregularity of $\infty$ is introduced according to the parabolic measure of $\infty$ is null or positive, and the necessary and sufficient condition for the Problem $\tilde{\mathcal{A}}_0$ is proved in terms of the Wiener-type criterion for the regularity of $\infty$.

In the probabilistic context, the formulated problems $\mathcal{A}_\gamma$ and $\tilde{\mathcal{A}}_\gamma$ are generalizations of the Kolmogorov problem to establish asymptotic laws for the *h-Brownian processes* [15].

1.1. **Overture: Kolmogorov-Petrovsky-type test for the Removability of the Fundamental Singularity.** For the special case of domains

(1.11) $$\Omega = \{(x,t) : |x - \gamma| < l(t), t > 0\} \subset \mathbb{R}_+^{N+1},$$

where $l \in C(\mathbb{R}_+; \mathbb{R}_+)$, $l(0) = 0$, $l(t) > 0$ for $t > 0$, the solution of the **Problem** $\mathcal{A}_\gamma$ reads:

**Theorem 1.2.** *$u_* \equiv 0$ or $u_* > 0$, that is to say, the fundamental singularity is removable or non-removable according to the following integral diverges or converges*

$$\int_{0+} t^{-\frac{N}{2}-1} l^{\frac{N}{2}}(t) e^{-\frac{l}{4t}}\, dt \tag{1.12}$$

The result is local. The removability of the fundamental singularity is dictated by the boundary of the domain near the singularity point. Precisely, it is defined by the thinness of the exterior set $\mathbb{R}^{N+1}_+ - \Omega$ near the singularity point $\mathcal{O}$.

The removability of the singularity is locally order-preserving. Precisely, if for some $\delta > 0$ we have $\Omega_1 \cap \{0 < t < \delta\} \subset \Omega_2 \cap \{0 < t < \delta\}$, then removability of the fundamental singularity for $\Omega_2$ (or non-removability for $\Omega_1$) implies the same for $\Omega_1$ (or $\Omega_2$) (Lemma 4.2, Section 4).

An equivalent form of the criterion can be written if we choose $l(t) = (4t \log \rho(t))^{\frac{1}{2}}$, and consider the domain $\Omega$ such that

$$\Omega \cap \{0 < t < \delta\} = \{(x,t): |x - \gamma|^2 < 4t \log \rho(t), 0 < t < \delta\} \tag{1.13}$$

such that

$$\rho \in C(0, \delta),\ \rho > 0;\ \rho(t) \uparrow +\infty,\ t \log \rho(t) \to 0,\ \text{as } t \downarrow 0. \tag{1.14}$$

Then the claim of the Theorem 1.2 remains valid if the integral (1.12) is replaced with

$$\int_{0+} \frac{|\log \rho(t)|^{\frac{N}{2}}}{t\rho(t)}\, dt \tag{1.15}$$

Some examples of functions $\rho$ with divergent integral (1.15) are as follows:

$$|\log t|,\ |\log t| \cdot \log_2^{\frac{N}{2}+1} t,\ |\log t| \cdot \log_2^{\frac{N}{2}+1} t \cdot \log_3 t \cdots \log_k t,\ k = 3, 4, \ldots \tag{1.16}$$

On the other side, for $\forall \epsilon > 0$, the integral (1.15) converges for the corresponding functions

$$|\log t|^{1+\epsilon},\ |\log t| \cdot \log_2^{\frac{N}{2}+1+\epsilon} t,\ |\log t| \cdot \log_2^{\frac{N}{2}+1} t \cdot \log_3 t \cdots \log_k^{1+\epsilon} t,\ k = 3, 4, \ldots \tag{1.17}$$

We adopt the notation

$$\log_2 t = \log|\log t|,\ \log_k t = \log \log_{k-1} t,\ k = 3, 4, \ldots$$

Hence, we have the following law for the removability of the fundamental singularity. For arbitrary integer $k \geq 4$, consider a domain

$$|x - \gamma|^2 \leq 4t \left[\log_2 t + \left(\frac{N}{2} + 1\right) \log_3 t + \cdots + \log_k^{1+\epsilon} t\right] \tag{1.18}$$

Then $u_* \equiv 0$ or $u_* > 0$, that is to say, the fundamental singularity is removable or non-removable according to $\epsilon = 0$ or $\epsilon > 0$.

**Probabilistic counterpart:** Let $\{x(t) = (x_1(t), \ldots, x_N(t)) : t \geq 0, P_\bullet\}$ be an $N$-dimensional $h$-Brownian process, and $P_\bullet(B)$ is a probability of the event $B$ as a function of the starting point $x(\tau)$ with $\tau > 0$ [15]. $h$-Brownian motion $x(t)$ from a point $x(\tau)$ is an almost surely continuous process whose sample functions never leave $\mathbb{R}^{N+1}_+$ and proceed downward, that is, in the direction of decreasing $t$. In fact, almost every path starting at $x(\tau)$ has a finite lifetime $\tau$ and tends to the boundary point $\mathcal{O}$ as $t \downarrow 0$ [15]. Consider a radial part $r(t) = |x(t) - \gamma| : 0 < t \leq \tau$ of the $h$-Brownian path starting at $x(\tau)$. Kolmogorov's 01 law implies that $P_{x(\tau)}[r(t) < h(t), t \downarrow 0] = 0$ or $1$; $h$ is said to belong the upper class if the probability is 1 and to the lower class otherwise. The probabilistic analog of Theorem 1.2 states that if $l \in \uparrow$ and if $t^{-1/2}l \in \downarrow$ for small $t > 0$, then $l$ belongs to the upper class or to the lower class according as the integral (1.12) converges or diverges. In particular, for arbitrary integer $k \geq 4$

$$l(t) = 2t^{\frac{1}{2}}\left[\log_2 t + \left(\frac{N}{2} + 1\right)\log_3 t + \cdots + \log_k^{1+\epsilon} t\right]^{\frac{1}{2}}$$

belongs to the upper or lower class according as $\epsilon > 0$ or $\epsilon \leq 0$.

Next, we describe the solution of the Problem $\tilde{\mathcal{A}}_\gamma$ for a special class of domains

(1.19) $\qquad \tilde{\Omega} = \{(x,t) : |x + 2t\gamma|^2 < -4t\log\rho(t), -\infty < t < \delta\} \subset \mathbb{R}_-^{N+1},$

where $\delta \ll -1$ and

(1.20) $\qquad \rho \in C(-\infty, \delta), \ \rho > 0, \ \rho(t) \to +\infty, \ t^{-1}\log\rho(t) \to 0, \text{ as } t \downarrow -\infty.$

**Theorem 1.3.** *$\tilde{u}_* \equiv 0$ or $\tilde{u}_* > 0$, that is to say the singularity at $\infty$ is removable or non-removable according to the following integral diverges or converges*

(1.21) $$\int_{-\infty} \frac{|\log\rho(t)|^{\frac{N}{2}}}{t\rho(t)} dt$$

Typical examples for the divergence or the convergence of the integral (1.21) are given by (1.16),(1.17) just by replacing $t$ with $|t|$. Hence, we have the following law for the removability of singularity at $\infty$. For arbitrary integer $k \geq 4$, consider a domain (1.19) with

$$|x + 2t\gamma|^2 < -4t\left[\log_2 |t| + \left(\frac{N}{2} + 1\right)\log_3 |t| + \cdots + \log_k^{1+\epsilon} |t|\right]$$

Then $u_* \equiv 0$ or $u_* > 0$, that is to say the singularity at $\infty$ is removable or non-removable according to $\epsilon \leq 0$ or $\epsilon > 0$.

Remarkably, in the particular case $\gamma = 0$, Theorem 1.3 coincides with the celebrated Kolmogorov-Petrovski test [27, 4].

**Probabilistic counterpart:** Let $\{x(t) = (x_1(t), ..., x_N(t)) : t < 0, P_\bullet\}$ be an $N$-dimensional $\tilde{h}$-Brownian process, and $P_\bullet(B)$ is a probability of the event $B$ as a function of the starting point $x(\tau)$ with $\tau < 0$ [15]. $\tilde{h}$-Brownian motion $x(t)$ from a point $x(\tau)$ is an almost surely continuous process whose sample functions never leave $\mathbb{R}_-^{N+1}$ and proceed downward, that is, in the direction of decreasing $t$. Almost every path starting at $x(\tau)$ tends to the boundary point $\infty$ as $t \downarrow -\infty$ [15]. Consider a radial part $r(t) = |x(t) + 2t\gamma| : t < 0$ of the $\tilde{h}$-Brownian path starting at $x(\tau)$. Kolmogorov's 01 law implies that $P_{x(\tau)}[r(t) < l(t), t \downarrow -\infty] = 0$ or $1$; $l$ is said to belong to the upper class if the probability is 1 and to the lower class otherwise. The probabilistic analog of Theorem 1.3 states that if $\rho$ satisfies (1.20), then $l(t) = 2(-t\log\rho(t))^{\frac{1}{2}}$ belongs to the upper class or the lower class according as the integral (1.21) converges or diverges. In particular, for arbitrary integer $k \geq 4$

$$l(t) = 2|t|^{\frac{1}{2}}\left[\log_2 |t| + \left(\frac{N}{2} + 1\right)\log_3 |t| + \cdots + \log_k^{1+\epsilon} |t|\right]^{\frac{1}{2}}$$

belongs to the upper or lower class according as $\epsilon > 0$ or $\epsilon \leq 0$.

## 2. Formulation of Problems

Being a generalization of the Kolmogorov problem, the Problems $\mathcal{A}_\gamma$, $\tilde{\mathcal{A}}_\gamma$, and their solution expressed in Theorems 1.2 and 1.3 has far-reaching measure-theoretical, topological and probabilistic implications in Analysis, PDEs and Potential theory. The goal of this section is to formulate three outstanding problems equivalent to the Problems $\mathcal{A}_\gamma$ and $\tilde{\mathcal{A}}_\gamma$. Since the problems $\mathcal{A}_\gamma$ and $\tilde{\mathcal{A}}_\gamma$ are equivalent via Appell transformation, without loss of generality we are going to formulate the problems in the framework of the Problem $\mathcal{A}_\gamma$. The equivalent formulation can be pursued in the framework of the Problem $\tilde{\mathcal{A}}_\gamma$ by replacing the triple $(\mathbb{R}_+^{N+1}, \Omega, h)$ with singularity point at $\mathcal{O}$ through the triple $(\mathbb{R}_-^{N+1}, \tilde{\Omega}, \tilde{h})$ with the singularity point at $\infty$ respectively.

### 2.1. Unique Solvability of the Singular Parabolic Dirichlet Problem.
Consider a *singular* **parabolic Dirichlet problem**(PDP) (1.2). The solution of the PDP can be constructed by Perron's method (or the method by Perron, Wiener, Brelot, and Bauer)[15, 30]. Let us introduce some necessary terminology.

We will often write a typical point $z \in \mathbb{R}^{N+1}$ as $z = (x,t), x \in \mathbb{R}^N, t \in \mathbb{R}$. A smooth solution of the heat equation is called a *parabolic function*. A bounded open set $U \subset \mathbb{R}^{N+1}$ is *regular* if for each continuous $f : \partial U \to \mathbb{R}$ there exists one (and only one) parabolic function $H_f^U : U \to \mathbb{R}$, such that

$$\lim_{\substack{z \to w \\ w \in \partial U}} H_f^U(z) = f(w), \ w \in \partial U.$$

A function $u$ is called a *superparabolic* in $\Omega$ if it satisfies the following conditions:

(1) $-\infty < u \leq +\infty$, $u < +\infty$ on a dense subset of $\Omega$;
(2) $u$ is lower semicontinuous (l.s.c.);
(3) for each regular open set $U \subset \Omega$ and each parabolic function $v \in C(\bar{U})$, the inequality $u \geq v$ on $\partial U$ implies $u \geq v$ in $U$.

A function $u$ is called a *subparabolic* if $-u$ is a superparabolic.

A function $u = v/h$ is called a *h-parabolic*, *h-superparabolic*, or *h-subparabolic* in $\Omega$ if $v$ is *parabolic*, *superparabolic*, or *subparabolic* [15].

We use the notation $\mathcal{S}_h(\Omega)$ for a class of all $h$-superparabolic functions in $\Omega$. Similarly, the class of all $h$-subparabolic functions in $\Omega$ is $-\mathcal{S}_h(\Omega)$.

Given boundary function $f$ on $\partial \Omega$, consider a $h$-**parabolic Dirichlet problem** ($h$-**PDP**): **find** $h$-**parabolic function** $u$ **in** $\Omega$ **such that**

$$(2.1) \qquad u = f \quad \text{on } \partial\Omega$$

It is easy to see that $h$-parabolic function $u = \frac{v}{h}$ is a bounded solution of the $h$-PDP if and only if $v$ is a solution of the PDP (1.2).

Assuming for a moment that $f \in C(\partial\Omega)$, the *generalized upper (or lower) solution* of the $h$-PDP is defined as

$$(2.2) \qquad {}^h\bar{H}_f^\Omega \equiv \inf\{u \in \mathcal{S}_h(\Omega) : \liminf_{z \to w, z \in \Omega} u \geq f(w) \quad \text{for all } w \in \partial\Omega\}$$

$$(2.3) \qquad {}^h\underline{H}_f^\Omega \equiv \sup\{u \in -\mathcal{S}_h(\Omega) : \limsup_{z \to w, z \in \Omega} u \leq f(w) \quad \text{for all } w \in \partial\Omega\}$$

The class of functions defined in (2.2) (or in (2.3)) is called *upper class* (or *lower class*) of the $h$-PDP. According to classical potential theory [15], $f$ is a $h$-*resolutive* *boundary function* in the sense that
$$^h\bar{H}_f^\Omega \equiv {^h\underline{H}_f^\Omega} \equiv {^hH_f^\Omega}.$$

The indicator function of any Borel measurable boundary subset, and equivalently any bounded Borel measurable boundary function is resolutive. Being $h$-parabolic in $\Omega$, $^hH_f^\Omega$ is called a generalized solution of the $h$-PDP for $f$. The generalized solution is unique by construction. It is essential to note that the construction of the generalized solution is accomplished by prescribing the behavior of the solution at $\mathcal{O}$.

Equivalently, we can define a generalized solution of the PDP (1.2):

**Definition 2.1.** Let $g : \partial\Omega \to \mathbb{R}$ be a boundary function, such that $g/h$ is a bounded Borel measurable. Then $g$ is called a resolutive boundary function for the PDP (1.2), if $f = g/h$ (extended to $\mathcal{O}$) is $h$-resolutive for the $h$-PDP. The function

$$(2.4) \qquad H_g^\Omega := h \, {^hH_f^\Omega}$$

is called a generalized solution of the PDP (1.2).

Again, note that the unique solution $H_g^\Omega$ of the PDP (1.2) is constructed by prescribing the behavior of the ratio $H_g^\Omega/h$ at $\mathcal{O}$.

The elegant theory, while identifying a class of unique solvability, leaves the following questions open:

- Would a unique solution of the $h$-DP still exist if its limit at $\mathcal{O}$ were not specified? That is, could it be that the solutions would pick up the "boundary" value $f(\mathcal{O})$ without being required? Equivalently, would a unique solution of the PDP (1.2) still exist if the limit of the ratio of solution to $h$ at $\mathcal{O}$ is not prescribed? In particular, is the fundamental singularity at $\mathcal{O}$ removable?
- What if the boundary datum $f$ (or $g/h$) on $\partial\Omega$, while being continuous at $\partial\Omega \setminus \{\mathcal{O}\}$, does not have a limit at $\mathcal{O}$, for example, it exhibits bounded oscillations. Is the $h$-PDP (or PDP (1.2)) uniquely solvable?

**Example 2.2.** Let $\Omega = \mathbb{R}_+^{N+1}$. It is easy to see that the boundary of $\mathbb{R}_+^{N+1}$ is $h$-resolutive and the only possible solutions of the $h$-PDP in $\mathbb{R}_+^{N+1}$ are constants. Precisely, the unique solution of the $h$-PDP is identical with the constant $f(\mathcal{O})$. Indeed, for arbitrary $\epsilon > 0$, the function

$$u(\cdot) = f(\mathcal{O}) + \frac{\epsilon}{h(\cdot)} \quad \left( \text{ or } v(\cdot) = f(\mathcal{O}) - \frac{\epsilon}{h(\cdot)} \right)$$

is in the upper class (or lower class) for $h$-PDP in $\mathbb{R}_+^{N+1}$ for $f$. Hence,

$$f(\mathcal{O}) - \frac{\epsilon}{h(\cdot)} \leq {^h\underline{H}_f^{\mathbb{R}_+^{N+1}}}(\cdot) \leq {^h\bar{H}_f^{\mathbb{R}_+^{N+1}}}(\cdot) \leq f(\mathcal{O}) + \frac{\epsilon}{h(\cdot)}$$

Since $\epsilon > 0$ is arbitrary, the assertion follows. Equivalently, all possible solutions of the PDP (1.2)in $\mathbb{R}_+^{N+1}$ are constant multiples of $h$, and the unique solution is identified by prescribing the ratio $u/h$ at $\mathcal{O}$.

Example 2.2 demonstrates that if $\Omega = \mathbb{R}_+^{N+1}$, the answer is negative and arbitrary constant $C$ is a solution of the $h$-PDP, $Ch$ is a solution of the PDP (1.2), and the fundamental singularity at $\mathcal{O}$ is not removable.

The positive answer to these fundamental questions is possible if $\Omega$ is not too sparse, or equivalently $\Omega^c \cap \mathbb{R}_+^{N+1}$ is not too thin near $\mathcal{O}$. The principal purpose of this paper is to prove the necessary and sufficient condition for the non-thinness of $\Omega^c \cap \mathbb{R}^{N+1}_+$ near $\mathcal{O}$ which is equivalent to the uniqueness of the solution of the $h$-PDP (or PDP ((1.2)) without specification of the boundary function (or ratio of the boundary function to $h$) at $\mathcal{O}$.

Furthermore, given *bounded Borel measurable* function $f = g/h : \partial\Omega \setminus \{\mathcal{O}\} \to \mathbb{R}$, we fix an arbitrary finite real number $\bar{f}$, and extend a function $f$ as $f(\mathcal{O}) = \bar{f}$. The extended function is a bounded Borel measurable on $\partial\Omega$ and there exists a unique solution ${}^h H^\Omega_f$ of the $h$-PDP, and the unique solution of the PDP (1.2) is given by (2.4). The major question now becomes:

**Problem 1:** *How many bounded solutions do we have, or does the constructed solution depend on $\bar{f}$ ?*

## 2.2. Characterization of the $h$-Parabolic Measure of Singularity Point.

For a given boundary Borel subset $A \subset \partial\Omega$, denote the indicator function of $A$ as $1_A$. Indicator functions of the Borel measurable subsets of $\partial\Omega$ are resolutive ([15]). *$h$-Parabolic measure* of the boundary Borel subset $A$ is defined as ([15]):

$$\mu^h_\Omega(z, A) = {}^h H^\Omega_{1_A}(z),$$

where $z \in \Omega$ is a reference point. It is said that $A$ is an $h$-parabolic measure null set if $\mu_\Omega(\cdot, A)$ vanishes identically in $\Omega$. If this is not the case, $A$ is a set of positive $h$-parabolic measure. In particular, the $h$-parabolic measure of $\{\mathcal{O}\}$ is well defined:

$$\mu^h_\Omega(\cdot, \{\mathcal{O}\}) = {}^h H^\Omega_{1_{\{\mathcal{O}\}}}(\cdot).$$

The following formula is true for the solution ${}^h H^\Omega_f$ of the $h$-PDP [15]:

$$(2.5) \qquad {}^h H^\Omega_f(z) = \int_{\partial\Omega} f(w) \mu^h_\Omega(z, dw), \quad z \in \Omega$$

Since $f$ is extended to $\{\mathcal{O}\}$ as $f(\mathcal{O}) = \bar{f}$, we have

$$(2.6) \qquad {}^h H^\Omega_f(z) = \int_{\partial\Omega \setminus \{\mathcal{O}\}} f(w) \mu^h_\Omega(z, dw) + \bar{f} \cdot {}^h H^\Omega_{1_{\{\mathcal{O}\}}}(z), \quad z \in \Omega$$

This formula implies that the uniqueness of the solution to the $h$-PDP without prescribing the behavior of the solution at the singularity point $\mathcal{O}$, that is to say, the independence of ${}^h H^\Omega_f$ on $\bar{f}$ is equivalent to whether or not $\mathcal{O}$ is an $h$-parabolic measure null set. Equivalently, according to the formula (2.4) the following formula is true for the unique solution of the PDP (1.2):

$$(2.7) \qquad H^\Omega_g(z) = h(z) \int_{\partial\Omega} \frac{g(w)}{h(w)} \mu^h_\Omega(z, dw), \quad z \in \Omega$$

Splitting the integral as in (2.6) we have

$$(2.8) \qquad H^\Omega_g(z) = h(z) \int_{\partial\Omega \setminus \{\mathcal{O}\}} \frac{g(w)}{h(w)} \mu^h_\Omega(z, dw) + \bar{f} h(z) \, {}^h H^\Omega_{1_{\{\mathcal{O}\}}}(z), \quad z \in \Omega.$$

where $\bar{f}$ is a prescribed limit value of $H^\Omega_g / h$ at $\mathcal{O}$. Similar to its counterpart (2.6), the formula (2.8) demonstrates that the uniqueness of the solution $u$ of the PDP (1.2) without prescribing $u/h$ at $\mathcal{O}$ is equivalent to whether or not $\mathcal{O}$ is an $h$-parabolic measure null set.

Hence, the following problem is the measure-theoretical counterpart of the Problem 1:

**Problem 2:** *Given $\Omega$, is the $h$-parabolic measure of $\{\mathcal{O}\}$ null or positive ?*

From the Example 2.2 demonstrated above it follows that in the particular case with $\Omega = \mathbb{R}_+^{N+1}$, we have

$$(2.9) \qquad \mu^h_{\mathbb{R}_+^{N+1}}(\cdot, \{\mathcal{O}\}) \equiv 1, \quad \mu^h_{\mathbb{R}_+^{N+1}}(\cdot, \partial\mathbb{R}_+^{N+1} - \{\mathcal{O}\}) \equiv 0.$$

**Example 2.3.** For arbitrary $c > 0$ consider a domain bounded by the level set of $h$:

$$\Omega = \{h > c'\} = \{z : |x - \gamma|^2 < -2Nt\log\frac{t}{c}, \ 0 < t < c\},$$

where $c' = (4\pi c)^{-\frac{N}{2}}$. It is easy to see that the $h$-harmonic measure of $\{\mathcal{O}\}$ is positive, and we have

$$(2.10) \qquad {}^hH^\Omega_{1_{\{\mathcal{O}\}}}(z) = \frac{h(z) - c'}{h(z)}, \ z \in \Omega.$$

Both **Problem 1 and 2 are equivalent to the Problem $\mathcal{A}_\gamma$ formulated in Section 1**. The connection follows from the following formula:

$$(2.11) \qquad u_*(z) = h(z) \ {}^hH^\Omega_{1_{\{\mathcal{O}\}}}(z), \ z \in \Omega$$

To establish (2.11), first note that the $h$-parabolic function $u_*/h$ is in the lower class of the Perron's solution ${}^hH^\Omega_{1_{\{\mathcal{O}\}}}$, which imply that

$$(2.12) \qquad \frac{u_*(z)}{h(z)} \leq \ {}^hH^\Omega_{1_{\{\mathcal{O}\}}}(z), \ z \in \Omega.$$

Moreover, ${}^hH^\Omega_{1_{\{\mathcal{O}\}}}$ itself is in the lower class of the Perron's solution $u_n/h$ of the $h$-PDP in $\Omega_n$ with boundary function $1_{\partial\Omega_n \cap \{t=n^{-1}\}}$, where $u_n$ is a solution of PDP (1.3). Therefore, we have

$$(2.13) \qquad {}^hH^\Omega_{1_{\{\mathcal{O}\}}}(z) \leq \frac{u_n(z)}{h(z)}, \ z \in \Omega_n.$$

passing to the limit as $n \to \infty$, from (1.5), (2.12) and (2.13), (2.11) follows.

In light of the measure-theoretical counterpart of the removability of the fundamental singularity, we introduce a concept of $h$-regularity of the boundary point $\mathcal{O}$.

**Definition 2.4.** $\mathcal{O}$ is said to be *$h$-regular* for $\Omega$ if it is an $h$-parabolic measure null set. Conversely, $\mathcal{O}$ is *$h$-irregular* if it has a positive $h$-parabolic measure.

Hence, Theorem 1.2 establishes a criterion for the removability of the fundamental singularity in terms of the necessary and sufficient condition for the *$h$-regularity* of $\mathcal{O}$.

Similarly, in the context of the singular PDP (1.7), and corresponding $\tilde{h}$-PDP, we have the formulae analogous to (2.6), (2.8) and (2.11):

$$\text{(2.14)} \quad {}^{\bar{h}}H_f^{\tilde{\Omega}}(z) = \int_{\partial \tilde{\Omega}\setminus\{\infty\}} f(w)\mu_{\tilde{\Omega}}^{\bar{h}}(z,dw) + \bar{f}\cdot {}^{\bar{h}}H_{1_{\{\infty\}}}^{\tilde{\Omega}}(z), \quad z \in \tilde{\Omega}$$

$$\text{(2.15)} \quad H_g^{\tilde{\Omega}}(z) = \tilde{h}(z) \int_{\partial \tilde{\Omega}\setminus\{\infty\}} \frac{g(w)}{h(w)}\mu_{\tilde{\Omega}}^{\bar{h}}(z,dw) + \bar{f}\tilde{h}(z)\, {}^{\bar{h}}H_{1_{\{\infty\}}}^{\tilde{\Omega}}(z), \quad z \in \tilde{\Omega}$$

$$\text{(2.16)} \quad \tilde{u}_*(z) = \tilde{h}(z)\, {}^{\bar{h}}H_{1_{\{\infty\}}}^{\tilde{\Omega}}(z), \ z \in \tilde{\Omega}.$$

We introduce a concept of $\tilde{h}$-regularity of the boundary point $\infty$ for $\tilde{\Omega} \subset \mathbb{R}_-^{N+1}$.

**Definition 2.5.** $\infty$ is said to be $\tilde{h}$-*regular* for $\tilde{\Omega}$ if it is an $\tilde{h}$-parabolic measure null set. Conversely, $\infty$ is $\tilde{h}$-*irregular* if it has a positive $\tilde{h}$-parabolic measure.

In fact, in the particular case with $\gamma = 0, \tilde{h} \equiv 1$, it coincides with the concept of regularity of $\infty$ introduced in [1]. Theorem 1.3 presents a criterion for the removability of the fundamental singularity at $\infty$ in terms of the necessary and sufficient condition for the $\tilde{h}$-*regularity* of $\infty$.

2.3. **Boundary Regularity in Singular Dirichlet Problem.** The notion of the $h$-regularity of $\mathcal{O}$ is, in particular, relates to the notion of continuity of the solution to the $h$-PDP at $\mathcal{O}$.

**Problem 3:** *Given $\Omega$, whether or not*

$$\text{(2.17)} \quad \liminf_{z\to\mathcal{O}, z\in\partial\Omega} f \leq \liminf_{z\to\mathcal{O}, z\in\Omega} {}^{h}H_f^{\Omega} \leq \limsup_{z\to\mathcal{O}, z\in\Omega} {}^{h}H_f^{\Omega} \leq \limsup_{z\to\mathcal{O}, z\in\partial\Omega} f,$$

for all bounded $f \in C(\partial\Omega \setminus \{\mathcal{O}\})$.

Note that if $f$ has a limit at $\mathcal{O}$, (2.17) simply means that the solution ${}^{h}H_f^{\Omega}$ is continuous at $\mathcal{O}$.

The equivalent problem in the context of the PDP (1.2) is the following:

**Problem 3′:** *Given $\Omega$, whether or not*

$$\text{(2.18)} \quad \liminf_{z\to\mathcal{O}, z\in\partial\Omega} \frac{g}{h} \leq \liminf_{z\to\mathcal{O}, z\in\Omega} \frac{H_g^{\Omega}}{h} \leq \limsup_{z\to\mathcal{O}, z\in\Omega} \frac{H_g^{\Omega}}{h} \leq \limsup_{z\to\mathcal{O}, z\in\partial\Omega} \frac{g}{h},$$

for all $g$ such that $\frac{g}{h} \in C(\partial\Omega \setminus \{\mathcal{O}\})$ and bounded.

In particular, if $g/h$ has a limit at $\mathcal{O}$, (2.18) means that the limit of the ratio $H_g^{\Omega}/h$ at $\mathcal{O}$ exists and equal to the limit of $g/h$.

Similarly, in the context of the singular PDP (1.7), and corresponding $\tilde{h}$-PDP, we can express the $\tilde{h}$-regularity of $\infty$ in terms of the regularity of the solution at $\infty$ by replacing (2.17), (2.18) with the conditions

$$\text{(2.19)} \quad \liminf_{z\to\infty, z\in\partial\tilde{\Omega}} f \leq \liminf_{z\to\infty, z\in\tilde{\Omega}} {}^{\bar{h}}H_f^{\tilde{\Omega}} \leq \limsup_{z\to\infty, z\in\tilde{\Omega}} {}^{\bar{h}}H_f^{\tilde{\Omega}} \leq \limsup_{z\to\infty, z\in\partial\tilde{\Omega}} f,$$

$$\forall \text{ bounded } f \in C(\partial\tilde{\Omega})$$

$$\text{(2.20)} \quad \liminf_{z\to\infty, z\in\partial\tilde{\Omega}} \frac{g}{\tilde{h}} \leq \liminf_{z\to\infty, z\in\tilde{\Omega}} \frac{H_g^{\tilde{\Omega}}}{\tilde{h}} \leq \limsup_{z\to\infty, z\in\tilde{\Omega}} \frac{H_g^{\tilde{\Omega}}}{\tilde{h}} \leq \limsup_{z\to\infty, z\in\partial\tilde{\Omega}} \frac{g}{\tilde{h}},$$

$$\forall \text{ bounded } \frac{g}{\tilde{h}} \in C(\partial\tilde{\Omega})$$

Theorem 1.2 (or 1.3) express the solutions to equivalent Problems 1-3 in terms of the Kolmogorov-Petrovsky-type criterion for the *h-regularity* of $\mathcal{O}$ (or *$\tilde{h}$-regularity* of $\infty$).

## 3. The Main Results

We now reformulate the main results of Theorems 1.2, 1.3 in a broader context as a solution of the equivalent Problems 1-3.

**Theorem 3.1.** *For arbitrary open set* $\Omega \subset \mathbb{R}^{N+1}_+$ *the following conditions are equivalent:*

(1) $\mathcal{O}$ *is h-regular (or h-irregular).*
(2) *Singular Parabolic Dirichlet Problem* (1.2), *and equivalently h-PDP has a unique (or infinitely many) solution(s).*
(3) *Boundary regularity conditions* (2.17), (2.18) *are satisfied (or, aren't satisfied).*
(4) *If* $\Omega$ *satisfies* (1.13),(1.14), *the integral* (1.15) *diverges (or converges).*

**Theorem 3.2.** *For arbitrary open set* $\tilde{\Omega} \subset \mathbb{R}^{N+1}_-$ *the following conditions are equivalent:*

(1) $\infty$ *is $\tilde{h}$-regular (or $\tilde{h}$-irregular).*
(2) *Singular Parabolic Dirichlet Problem* (1.7), *and equivalently $\tilde{h}$-PDP has a unique (or infinitely many) solution(s).*
(3) *Boundary regularity conditions* (2.19), (2.20) *are satisfied (or, aren't satisfied).*
(4) *If* $\tilde{\Omega}$ *satisfies* (1.19),(1.20), *the integral* (1.21) *diverges (or converges).*

3.1. **Historical Comments.** The major problem in the Analysis of PDEs is understanding the nature of singularities of solutions to the PDEs reflecting the natural phenomena. It would be convenient to make some remarks on the analysis of singularities for the Laplace and heat equations, as well as more general second-order elliptic and parabolic PDEs. The solvability, in some generalized sense, of the classical DP in a bounded open set $E \subset \mathbb{R}^N$, with prescribed data on $\partial E$, is realized within the class of resolutive boundary functions, identified by Perron's method and its Wiener [31, 32] and Brelot [13] refinements. Such a method is referred to as the PWB method, and the corresponding solutions are PWB solutions. Paralleling the theory of PWB solutions, the DP for the heat equation in an arbitrary open set is solvable within the class of resolutive boundary functions. We refer to [30, 15] for an account of the theory. Wiener, in his pioneering works [31, 32], proved a necessary and sufficient condition for the finite boundary point $x_o \in \partial E$ to be regular in terms of the "thinness" of the complementary set in the neighborhood of $x_o$. If the boundary of the domain is a graph in a neighborhood of $x_0$, the Wiener criterion is entirely geometrical. A key advance made in Wiener's work was an introduction of the concept of capacity - sub-additive set function dictated by the Laplacian for the accurate measuring of the thinness of the complementary set in the neighborhood of $x_0$ for the boundary regularity of harmonic function. Formalized through the powerful Choquet capacitability theorem [14], the concept of capacity became a standard tool for the characterization of singularities for the elliptic and parabolic equations. The question of removability of isolated singularities for the linear second-order elliptic and parabolic PDEs was settled in [28], and in [8, 9, 10]. Wiener criterion for the boundary continuity of harmonic functions became a canonical result driving the boundary regularity theory for the elliptic and parabolic PDEs. In 1935, Petrovsky proved a geometric necessary and sufficient condition for the regularity of the characteristic top boundary point for the heat equation in the domain of revolution [27] (see also [3]). In the same paper, he also presented an elegant solution of the Kolmogorov problem (see Section 1, Problem $\tilde{\mathcal{A}}_\gamma$) for the special domain of revolution (see also [4]). The results formed the so-called Kolmogorov-Petrovsky test for the asymptotic behavior of the standard Brownian path as $t \downarrow 0$ and $t \uparrow +\infty$, and opened a path for the deep connection between the regularity theory of elliptic and parabolic PDEs and asymptotic properties of the associated Markov processes [22]. The geometric iterated logarithm test for the regularity of the boundary point for an arbitrary open set with respect to heat equation is proved in [2]. Paralleling the Wiener regularity theory, Wiener's criterion for the regularity of the finite boundary point for the heat equation was formulated in [24] along with the proof of the irregularity assertion. The problem was accomplished in [16], where the long-awaited regularity assertion was proved. As in its elliptic counterpart, the concept of *heat capacity* was a key concept to extend the Wiener regularity theory to the case of heat equation [30]. However, the major technical difficulty in doing so was connected to the nature of singularities of the fundamental solution of the Laplace and heat equations. The former is an isolated singularity for the spherical level sets of the fundamental solution, while the latter is a non-isolated singularity point for the level sets of the fundamental solution of the heat equation. To complete the Wiener regularity result at finite boundary points for the heat equation, the major technical advance of paper [16] was a proof of elegant boundary Harnack estimate near the non-isolated singularity point of the level sets of the fundamental solution to the backward heat equation. The result of [16] was extended to the class of linear second-order divergence form parabolic PDEs with $C^1$-Dini continuous coefficients in [20, 18].

In [25] it is proved that the Wiener test for the regularity of finite boundary points concerning second-order divergence form uniformly elliptic operator with bounded measurable coefficients coincides with the classical Wiener test for the boundary regularity of harmonic functions. The Wiener test for the regularity of finite boundary points for linear degenerate elliptic equations is proved in [17]. The Wiener test for the regularity of finite boundary points for quasilinear elliptic equations was settled due to [26, 19, 23]. Nonlinear potential theory was developed along the same lines as classical potential theory for the Laplace operator, for which we refer to [21].

To solve the DP in an unbounded open set, Brelot introduced the idea of compactifying $\mathbb{R}^N$ into $\mathbb{R}^N \cup \{\infty\}$, where $\infty$ is the point at $\infty$ of $\mathbb{R}^N$ [12]. PWB-method is extended to the compactified framework, thus providing a powerful existence and uniqueness result for the DP in arbitrary open sets in the class of resolutive boundary functions. The new concept of regularity of $\infty$ was introduced in [5] for the classical DP, and in [6] for its parabolic counterpart. The DP with bounded Borel measurable boundary function has one and only one or infinitely many solutions without prescribing the boundary value at $\infty$. The point at $\infty$ is called a regular if there is a unique solution, and it is called irregular otherwise. Equivalently, in the measure-theoretical context, the new concept of regularity or irregularity of $\infty$ is introduced according to whether the harmonic measure of $\infty$ is null or positive. In [5] the Wiener criterion for the regularity of $\infty$ in the classical DP for the Laplace equation in an open set $E \subset \mathbb{R}^N$ with $N \geq 3$ is proved. In [6] it is proved that the Wiener criterion at $\infty$ for the linear second-order divergence form elliptic PDEs with bounded measurable coefficients coincides with the Wiener criterion at $\infty$ for the Laplacian operator. The Wiener criterion at $\infty$ for the heat equation is proved in [1]. Remarkably, the Kolmogorov problem (see Section 1, Problem $\tilde{\mathcal{A}}_\gamma$) is a particular case of the problem of uniqueness of the bounded solution of the parabolic Dirichlet problem in arbitrary open set in $\mathbb{R}^{N+1}$ without prescribing the limit of the solution at $\infty$. Hence, the Wiener criterion at $\infty$ proved in [1] presents a full solution to the Kolmogorov problem.

The new Wiener criterion at $\infty$ for the elliptic and parabolic PDEs broadly extends the role of the Wiener regularity theory in a classical Analysis. The Wiener test at $\infty$ arises as a global characterization of uniqueness in boundary value problems in arbitrary unbounded open sets. From a topological point of view, the Wiener test at $\infty$ arises as a thinness criterion at $\infty$ in fine topology. In a probabilistic context, the Wiener test at $\infty$ characterizes asymptotic laws for the Markov processes whose generator is a given differential operator. The counterpart of the new Wiener test at a finite boundary point leads to uniqueness in a Dirichlet problem for a class of unbounded functions growing at a certain rate near the boundary point; a criterion for the removability of singularities and/or for unique continuation at the finite boundary point: let $E \subset \mathbb{R}^N, N \geq 3$ be an open set, and $x_0 \in E$ be a finite boundary point. Consider a singular Dirichlet problem for the linear second order uniformly elliptic PDE with bounded measurable coefficients in a class $O(|x - x_0|^{2-N}$ as $x \to x_0$. In [6] it is proved that the Wiener test at $x_0$ is a necessary and sufficient condition for the unique solvability of the singular Dirichlet problem, and equivalently for the removability of the fundamental singularity at $x_0$. In a recent paper [7] an appropriate 2D analog of this result is established. Let $E \subset \mathbb{R}^2$ be a Greenian open set, and $x_0 \in \partial E$ be a boundary point (finite or $\infty$). Consider a singular Dirichlet problem for the linear second-order uniformly elliptic operator with bounded measurable coefficients in the class $O(\log|x - x_0|)$ if $x_0$ is finite, and in a class of functions with logarithmic growth, if $x_0 = \infty$. In [7] it is proved that the Wiener criterion at $x_0$ is a necessary and sufficient condition for the unique solvability of the singular Dirichlet problem, and equivalently for the removability of the logarithmic singularity. Precisely, in [7] the concept of log-regularity (or log-irregularity) of the boundary point (finite or $\infty$) is introduced according as if log-harmonic measure of it is null or positive, and the removability of the logarithmic singularity is expressed in terms of the Wiener criterion for the log-regularity of $x_0$.

The goal of this paper is to establish a necessary and sufficient condition for the removability of the fundamental singularity, and equivalently for the unique solvability of the singular PDP. In this paper, we prove the Kolmogorov-Petrovsky-type test. We address the proof of the Wiener-type criterion in the forthcoming paper.

## 4. Preliminary Results

The equivalence of two problems formulated in $\mathbb{R}^{N+1}_+$ and $\mathbb{R}^{N+1}_-$ is a consequence of the **Appell transformation**. Consider a homeomorphism $A : \mathbb{R}^{N+1}_+ \cup \{\mathcal{O}\} \mapsto$

$\mathbb{R}_-^{N+1} \cup \{\infty\}$ with

(4.1) $\begin{cases} (x,t) \in \mathbb{R}_+^{N+1} \mapsto A(x,t) = \left(\frac{x}{2t}, -\frac{1}{4t}\right) \in \mathbb{R}_-^{N+1}; \ A(\mathcal{O}) = \infty \\ (x,t) \in \mathbb{R}_-^{N+1} \mapsto A^{-1}(x,t) = \left(-\frac{x}{2t}, -\frac{1}{4t}\right) \in \mathbb{R}_+^{N+1}; \ A^{-1}(\infty) = \mathcal{O} \end{cases}$

Let $\mathcal{P}(\Omega)$ be a class of parabolic functions in an open set $\Omega$. Given open set $\Omega \subset \mathbb{R}_+^{N+1}$, the Appell transformation is a homeomorphism $\mathbf{A} : \mathcal{P}(\Omega) \mapsto \mathcal{P}(A\Omega)$ defined as

(4.2) $\begin{cases} \Omega \subset \mathbb{R}_+^{N+1} : u \in \mathcal{P}(\Omega) \mapsto \mathbf{A}u(z) = (-\frac{\pi}{t})^{\frac{N}{2}} e^{-\frac{|x|^2}{4t}} u(A^{-1}(z)) \in \mathcal{P}(A\Omega) \\ \Omega \subset \mathbb{R}_-^{N+1} : v \in \mathcal{P}(\Omega) \mapsto \mathbf{A}^{-1}v(z) = F(z)v(A(z)) \in \mathcal{P}(A^{-1}\Omega) \end{cases}$

The claim follows from the following formula:

(4.3) $\begin{cases} \mathcal{H}[\mathbf{A}u(z)] = \frac{\pi^{N/2}}{4}(-t)^{-\frac{N}{2}-2} e^{-\frac{|x|^2}{4t}} \mathcal{H}[u(A^{-1}(z)], \ z \in A\Omega \subset \mathbb{R}_-^{N+1} \\ \mathcal{H}[\mathbf{A}^{-1}v(z)] = \frac{1}{4t^2} F(z)\mathcal{H}[v(A(z)], \ z \in A^{-1}\Omega \subset \mathbb{R}_+^{N+1}. \end{cases}$

In particular, the Appell transform of $h$ is given by

(4.4) $$\tilde{h}(x,t) = \mathbf{A}h(x,t) = e^{\langle x, \gamma \rangle + |\gamma|^2 t},$$

as it is defined in (1.6).

From the formula (4.3) it follows that The Appell transformation is a homeomorphism between $\mathcal{S}(\Omega) \cap C^2(\Omega)$ and $\mathcal{S}(A\Omega) \cap C^2(A\Omega)$, where $\mathcal{S}(\Omega)$ denotes the class of all superparabolic functions in $\Omega$. A simple approximation argument can be used to demonstrate that the hypothesis $C^2(\Omega)$ can be removed [15].

Appell transformation presents one-to-one mapping between the singular PDPs (1.2) and (1.7).

**Lemma 4.1.** (1) *Function $u$ is $h$-parabolic (or $h$-superparabolic) in open set $\Omega \subset \mathbb{R}_+^{N+1}$ if and only if $u(A^{-1}(x,t))$ is $\tilde{h}$-parabolic (or $\tilde{h}$-superparabolic) in $A\Omega \subset \mathbb{R}_-^{N-1}$.*

(2) $^hH_f^\Omega$ *is a solution of the $h$-parabolic Dirichlet problem in $\Omega \in \mathbb{R}_+^{N+1}$ if and only if $^hH_f^\Omega(A^{-1}(\cdot))$ is a solution of the $\tilde{h}$-parabolic Dirichlet problem in $A\Omega \subset \mathbb{R}_-^{N+1}$ with boundary function $f(A^{-1})$, i.e.*

(4.5) $$^hH_f^\Omega(A^{-1}(z)) = {}^{\tilde{h}}H_{f(A^{-1})}^{A\Omega}(z), \ z = (x,t) \in A\Omega.$$

(3) $\mathcal{O}$ *is $h$-regular for $\Omega \subset \mathbb{R}_+^{N+1}$ if and only if $\infty$ is $\tilde{h}$-regular for $A\Omega \subset \mathbb{R}_-^{N+1}$.*

(4) $H_g^\Omega$ *is a solution of the singular PDP (1.2) if and only if its Appell transform is a solution of the singular PDP (1.7) in $A\Omega$ with boundary function $\mathbf{A}g(A^{-1})$, i.e.*

(4.6) $$\mathbf{A}H_g^\Omega(z) = H_{\mathbf{A}g(A^{-1})}^{A\Omega}(z), \ z \in A\Omega$$

(4.7) $$\mathbf{A}^{-1}H_{\tilde{g}}^{\tilde{\Omega}}(z) = H_{\mathbf{A}^{-1}\tilde{g}(A)}^{A^{-1}\Omega}(z), \ z \in A^{-1}\tilde{\Omega}$$

(5) *Problems $\mathcal{A}_\gamma|_\Omega$ and $\tilde{\mathcal{A}}_\gamma|_{A\Omega}$ are equivalent, i.e. $u_* \equiv 0$ if and only if $\tilde{u}_* \equiv 0$.*

*Proof.*

(1) Let $\Omega \subset \mathbb{R}^{N+1}_+$ be an open set, and $u$ be $h$-parabolic function on $\Omega$, i.e.
$$u(x,t) = \frac{v(x,t)}{h(x,t)}, \ (x,t) \in \Omega,$$
where $v$ is a parabolic function in $\Omega$. Considering the Appel transform of $v = uh$ we have
$$\mathbf{A}(uh)(x,t) = \left(-\frac{\pi}{t}\right)^{\frac{N}{2}} e^{-\frac{|x|^2}{4t}} u(A^{-1}(x,t))h(A^{-1}(x,t))$$
$$= e^{\frac{|x+2\gamma t|^2}{4t} - \frac{|x|^2}{4t}} u(A^{-1}(x,t)) = \tilde{h}(x,t)u(A^{-1}(x,t)), \ (x,t) \in A\Omega \subset \mathbb{R}^{N+1}_-,$$
which implies that $u(A^{-1}(x,t))$ is $\tilde{h}$-parabolic function in $A\Omega$. On the other side, let $\Omega \subset \mathbb{R}^{N+1}_-$ be an open set, and $u$ be $\tilde{h}$-parabolic function on $\Omega$, i.e.
$$u(x,t) = \frac{v(x,t)}{\tilde{h}(x,t)}, \ (x,t) \in \Omega,$$
where $v$ is a parabolic function in $\Omega$. Considering the inverse Appel transform of $v = u\tilde{h}$ we have
$$\mathbf{A}^{-1}(u\tilde{h})(x,t) = F(x,t)u(A(x,t))\tilde{h}(A(x,t))$$
$$= (4\pi t)^{-\frac{N}{2}} e^{-\frac{|x|^2}{4t}} e^{\frac{2\langle x,\gamma\rangle - |\gamma|^2}{4t}} u(A(x,t))$$
$$= h(x,t)u(A(x,t)), \ (x,t) \in A^{-1}\Omega \subset \mathbb{R}^{N+1}_+,$$
which implies that $u(A(x,t))$ is $h$-parabolic function in $A^{-1}\Omega$. The presented proof applies to smooth superparabolic functions without any changes. Using the standard smoothing, the proof is extended to $h$- and $\tilde{h}$-superparabolic functions as well.

(2) According to (ii) $^h H^\Omega_f(A^{-1}(x,t))$ is $\tilde{h}$-parabolic in $A\Omega$, and we only need to verify that the relations (2.2), (2.3) corresponding to $\tilde{h}$-parabolic Dirichlet problem are satisfied. For arbitrary $z = (x,t) \in A\Omega$, we have

$$^h \bar{H}^\Omega_f(A^{-1}(z)) \equiv \inf\{u \in \mathcal{S}_h(\Omega) : \liminf_{\substack{A^{-1}(z) \to w \\ A^{-1}(z) \in \Omega}} u(A^{-1}(z)) \geq f(w) \ \forall w \in \partial\Omega\} \equiv$$

$$\inf\{u \in \mathcal{S}_{\tilde{h}}(A\Omega) : \liminf_{\substack{z \to w \\ z \in A\Omega}} u(z) \geq f(A^{-1}(w)) \ \forall w \in \partial A\Omega\} \equiv {}^{\tilde{h}} \bar{H}^{A\Omega}_{f(A^{-1})}(z),$$

and

$$^h \underline{H}^\Omega_f(A^{-1}(z)) \equiv \sup\{u \in -\mathcal{S}_h(\Omega) : \limsup_{\substack{A^{-1}(z) \to w \\ A^{-1}(z) \in \Omega}} u(A^{-1}(z)) \leq f(w) \ \forall w \in \partial\Omega\} \equiv$$

$$\sup\{u \in -\mathcal{S}_{\tilde{h}}(A\Omega) : \limsup_{\substack{z \to w \\ z \in A\Omega}} u(z) \leq f(A^{-1}(w)) \ \forall w \in \partial A\Omega\} \equiv {}^{\tilde{h}} \underline{H}^{A\Omega}_{f(A^{-1})}(z),$$

which implies (4.5).

(3) If $f = 1_{\{\mathcal{O}\}}$, (4.5) implies

(4.8) $$^h H^\Omega_{1_{\{\mathcal{O}\}}}(A^{-1}(z)) = {}^{\tilde{h}} H^{A\Omega}_{1_{\{\infty\}}}(z), \ z = (x,t) \in A\Omega,$$

which proves the claim.

(4) By using (2.4) and (4.5) we have

$$\mathbf{A}H_g^\Omega(z) = \mathbf{A}[h\ {}^h H_f^\Omega](z) = \Big(-\frac{\pi}{t}\Big)^{\frac{N}{2}} e^{-\frac{|x|^2}{4t}} h(A^{-1}(z))^h H_f^\Omega(A^{-1}(z))$$

(4.9) $$= \tilde{h}(z)\ {}^{\tilde{h}} H_{f(A^{-1})}^{A\Omega}(z) = H_{\tilde{h}f(A^{-1})}^{A\Omega}(z) = H_{\mathbf{A}g(A^{-1})}^{A\Omega}(z),\ z \in A\Omega$$

since

$$\mathbf{A}g(A^{-1}(z)) = \Big(-\frac{\pi}{t}\Big)^{\frac{N}{2}} e^{-\frac{|x|^2}{4t}} g(A^{-1}(z))$$

(4.10) $$= \Big(-\frac{\pi}{t}\Big)^{\frac{N}{2}} e^{-\frac{|x|^2}{4t}} h(A^{-1}(z))f(A^{-1}(z)) = \tilde{h}(z)f(A^{-1}(z)),\ z \in A\Omega.$$

Proof of the symmetric relation (4.7) is similar.

(5) The claim is a direct consequence of (2.11),(4.8). $\square$

The next lemma expresses the fact that the property of $h$-regularity of the singularity point is local and order-preserving.

**Lemma 4.2.** (1) If $\Omega_1 \subset \Omega_2 \subset \mathbb{R}_+^{N+1}$, then
  (a) ${}^h H_{1_{\{\mathcal{O}\}}}^{\Omega_2} \equiv 0 \implies {}^h H_{1_{\{\mathcal{O}\}}}^{\Omega_1} \equiv 0$;
  (b) ${}^h H_{1_{\{\mathcal{O}\}}}^{\Omega_1} \not\equiv 0 \implies {}^h H_{1_{\{\mathcal{O}\}}}^{\Omega_2} \not\equiv 0$;
(2) Let $\Omega \subset \mathbb{R}_+^{N+1}$, and $\Omega_\delta := \Omega \cap \{t < \delta\}, \delta > 0$. Then ${}^h H_{1_{\{\mathcal{O}\}}}^\Omega \equiv 0$ if and only if ${}^h H_{1_{\{\mathcal{O}\}}}^{\Omega_\delta} \equiv 0$ for some (and equivalently for all) $\delta > 0$.
(3) If $\tilde{\Omega}_1 \subset \tilde{\Omega}_2 \subset \mathbb{R}_-^{N+1}$, then
  (a) ${}^{\tilde{h}} H_{1_{\{\infty\}}}^{\tilde{\Omega}_2} \equiv 0 \implies {}^{\tilde{h}} H_{1_{\{\infty\}}}^{\tilde{\Omega}_1} \equiv 0$;
  (b) ${}^{\tilde{h}} H_{1_{\{\infty\}}}^{\tilde{\Omega}_1} \not\equiv 0 \implies {}^{\tilde{h}} H_{1_{\{\infty\}}}^{\tilde{\Omega}_2} \not\equiv 0$;
(4) Let $\tilde{\Omega} \subset \mathbb{R}_-^{N+1}$, and $\tilde{\Omega}_\delta := \tilde{\Omega} \cap \{t < -\delta\}, \delta > 0$. Then ${}^h H_{1_{\{\infty\}}}^{\tilde{\Omega}} \equiv 0$ if and only if ${}^h H_{1_{\{\infty\}}}^{\tilde{\Omega}_\delta} \equiv 0$ for some (and equivalently for all) $\delta > 0$.
(5) ${}^{\tilde{h}} H_{1_{\{\infty\}}}^{D_\delta} \equiv 0$ for $D_\delta = \{x \in \mathbb{R}^N, t > -\delta\}, \delta > 0$.

*Proof.*

(1) It is easy to see that ${}^h H_{1_{\{\mathcal{O}\}}}^{\Omega_2} = {}^h H_g^{\Omega_1}$ on $\Omega_1$, where

(4.11) $$g = \begin{cases} 1_{\{\mathcal{O}\}}, & \text{on } \partial\Omega_1 \cap \partial\Omega_2 \\ {}^h H_{1_{\{\mathcal{O}\}}}^{\Omega_2}, & \text{on } \partial\Omega_1 \cap \Omega_2. \end{cases}$$

Since Perron's solution is order-preserving, it follows that

(4.12) $$\ ^h H_{1_{\{\mathcal{O}\}}}^{\Omega_1} \leq {}^h H_g^{\Omega_1} \leq {}^h H_{1_{\{\mathcal{O}\}}}^{\Omega_2}, \text{ on } \Omega_1,$$

which implies the claims (1a) and (1b).

(2) The "only if" claim is trivial. To demonstrate the "if" claim, note that since $\partial\Omega_\delta \cap \{t = \delta\}$ is a parabolic measure null set for $\Omega_\delta$, we have

$$\ ^h H_{1_{\{\mathcal{O}\}}}^{\Omega}|_{\Omega_\delta} = {}^h H_{1_{\{\mathcal{O}\}}}^{\Omega_\delta} \equiv 0.$$

Therefore, by the maximum principle we have

$$\ ^h H_{1_{\{\mathcal{O}\}}}^{\Omega}|_{\Omega \setminus \overline{\Omega}_\delta} \equiv {}^h H_0^{\Omega \setminus \overline{\Omega}_\delta} \equiv 0,$$

which proves the claim (2).

The proof of (3) and (4) is identical to the proof of (1) and (2). The claim (5) is a consequence of the uniqueness result for the Cauchy problem [29].

## 5. Proofs of Main Results

**5.1. Proof of Theorem 1.2.** We assume that $\Omega$ satisfies (1.13), (1.14). Assume that the integral (1.15) converges. We aim to demonstrate that $u_* > 0$, or equivalently, $\{\mathcal{O}\}$ is $h$-irregular. Without loss of generality, we can assume that $\rho \in C^1(0, \delta)$. Indeed, otherwise we can select the function $\rho_1 \in C^1(0, \delta)$ which satisfy

$$\text{(5.1)} \qquad \frac{\rho(t)}{2} < \rho_1(t) < \rho(t), \ 0 < t < \delta,$$

and consider the domain

$$\text{(5.2)} \qquad \Omega_1 = \{(x,t): \ |x-\gamma|^2 < 4t \log \rho_1(t), 0 < t < \delta\}.$$

From (5.1) it follows that the integral (1.15) is convergent for $\rho_1$, and $\Omega_1 \subset \Omega \cap \{0 < t < \delta\}$. Therefore, $h$-irregularity of $\{\mathcal{O}\}$ for $\Omega_1$ would imply so for $\Omega$.

Consider a function
$$u(x,t) = 1 - \rho^{-1}(t) e^{\frac{|x-\gamma|^2}{4t}}$$
which is positive in $\Omega$ for $0 < \delta << 1$, vanishes on $\partial\Omega \cap \{|x-\gamma|^2 = 4t \log \rho(t)\}$, and satisfies

$$\text{(5.3)} \qquad \lim_{t \downarrow 0} u(\gamma, t) = \limsup_{(x,t) \to \mathcal{O}, (x,t) \in \Omega} u(x,t) = 1.$$

For a function $v = uh$ we have

$$\mathcal{H}v = h\left[\frac{\rho'(t)}{\rho^2(t)} e^{\frac{|x-\gamma|^2}{4t}} + \frac{N}{2} \frac{1}{t\rho(t)} e^{\frac{|x-\gamma|^2}{4t}} + \frac{1}{\rho(t)} \frac{|x-\gamma|^2}{2t^2} e^{\frac{|x-\gamma|^2}{4t}}\right]$$

$$\text{(5.4)} \qquad -(4\pi t)^{-\frac{N}{2}} \frac{1}{\rho(t)} \frac{|x-\gamma|^2}{2t^2} = \frac{\rho'(t)}{\rho^2(t)} (4\pi t)^{-\frac{N}{2}} + \frac{N}{2} \frac{1}{t\rho(t)} (4\pi t)^{-\frac{N}{2}}.$$

Now we construct a function $w$ with the following properties:

$$\text{(5.5)} \qquad \mathcal{H}w = -\frac{N}{2} \frac{1}{t\rho(t)} (4\pi t)^{-\frac{N}{2}}, \ w(x,t) < 0, \ \text{in } \Omega_n,$$

$$\text{(5.6)} \qquad \frac{w(\gamma, t)}{h(\gamma, t)} \geq -\frac{1}{2}, \ \text{for } n^{-1} \leq t \leq T,$$

for some fixed $0 < T < \delta$ and for all sufficiently large $n$. From (5.4),(5.5) it follows that the function $\tilde{u} = \frac{v+w}{h}$ is $h$-subparabolic, and by the maximum principle we have

$$\text{(5.7)} \qquad u_n(x,t) \geq \tilde{u}(x,t) \text{ in } \Omega_n$$

If we select the fixed value $T$ sufficiently small, it follows that $\tilde{u}(\gamma, T) > 1/3$. Therefore, within $\Omega_T = \Omega \cap \{t > T\}$, $\tilde{u}$ is greater than the function which is $h$-parabolic in $\Omega_T$, takes the value $1/4$ in $\{(x,t): |x-\gamma| \leq \epsilon, t = T\}$ for some $0 < \epsilon << 1$, and vanishes on the rest of the parabolic boundary of $\Omega_T$. Hence, we have a positive lower bound for $u_n$ in $\Omega_T$ which is independent of $n$. Obviously, the same lower bound holds for the limit function $u_*$, that is to say the $h$-parabolic measure of $\mathcal{O}$ would be positive. Thus, to complete the proof we need to construct the function $w$ with the properties (5.5),(5.6).

As a function $w$ we select a particular solution of the equation from (5.5):

$$(5.8) \quad w(x,t) = -\frac{N}{2} \int_{n^{-1}}^{t} \int_{|\xi-\gamma| \leq (4\tau \log \rho(\tau))^{\frac{1}{2}}} \frac{exp\left(-\frac{|x-\xi|^2}{4(t-\tau)}\right)}{(4\pi(t-\tau))^{\frac{N}{2}}} \frac{(4\pi\tau)^{-\frac{N}{2}}}{\tau \rho(\tau)} d\xi \, d\tau$$

We only need to check that (5.6) is satisfied. We have

$$\left|\frac{w(\gamma,t)}{h(\gamma,t)}\right| = \frac{N}{2} \int_{n^{-1}}^{t} \int_{|\xi-\gamma| \leq (4\tau \log \rho(\tau))^{\frac{1}{2}}} \frac{exp\left(-\frac{|\xi-\gamma|^2}{4(t-\tau)}\right)}{(4\pi(t-\tau))^{\frac{N}{2}}} \frac{1}{\tau \rho(\tau)} \left(\frac{\tau}{t}\right)^{-\frac{N}{2}} d\xi \, d\tau$$

We split the integral $\int_{n^{-1}}^{t} = \int_{t/2}^{t} + \int_{n^{-1}}^{t/2}$, and estimate the first one as follows:

$$\frac{N}{2} \int_{t/2}^{t} \int_{|\xi-\gamma| \leq (4\tau \log \rho(\tau))^{\frac{1}{2}}} \frac{exp\left(-\frac{|\xi-\gamma|^2}{4(t-\tau)}\right)}{(4\pi(t-\tau))^{\frac{N}{2}}} \frac{1}{\tau \rho(\tau)} \left(\frac{\tau}{t}\right)^{-\frac{N}{2}} d\xi \, d\tau \leq$$

$$N 2^{\frac{N}{2}-1} \int_{t/2}^{t} \frac{1}{\tau \rho(\tau)} \int_{|\xi-\gamma| \leq (4\tau \log \rho(\tau))^{\frac{1}{2}}} \frac{exp\left(-\frac{|\xi-\gamma|^2}{4(t-\tau)}\right)}{(4\pi(t-\tau))^{\frac{N}{2}}} d\xi \, d\tau \leq$$

$$(5.9) \quad N 2^{\frac{N}{2}-1} \pi^{-\frac{N}{2}} \int_{\mathbb{R}^N} e^{-|\xi|^2} d\xi \int_{t/2}^{t} \frac{d\tau}{\tau \rho(\tau)} = N 2^{\frac{N}{2}-1} \int_{t/2}^{t} \frac{d\tau}{\tau \rho(\tau)}.$$

From the convergence of the integral (1.15) it follows that the right-hand side of (5.9) converges to zero as $t \downarrow 0$. We then have

$$\frac{N}{2} \int_{n^{-1}}^{t/2} \int_{|\xi-\gamma| \leq (4\tau \log \rho(\tau))^{\frac{1}{2}}} \frac{exp\left(-\frac{|\xi-\gamma|^2}{4(t-\tau)}\right)}{(4\pi(t-\tau))^{\frac{N}{2}}} \frac{1}{\tau \rho(\tau)} \left(\frac{\tau}{t}\right)^{-\frac{N}{2}} d\xi \, d\tau \leq$$

$$\frac{N}{2} \int_{n^{-1}}^{t/2} \frac{1}{\tau \rho(\tau)} \left(\frac{\tau}{t}\right)^{-\frac{N}{2}} \int_{|\xi-\gamma| \leq (4\tau \log \rho(\tau))^{\frac{1}{2}}} \frac{exp\left(-\frac{|\xi-\gamma|^2}{4(t-\tau)}\right)}{(4\pi(t-\tau))^{\frac{N}{2}}} d\xi \, d\tau \leq$$

$$\omega_N \frac{N}{2} \int_{n^{-1}}^{t/2} \frac{1}{\tau \rho(\tau)} \left(\frac{\tau}{t}\right)^{-\frac{N}{2}} (2\pi t)^{-\frac{N}{2}} (4\tau \log \rho(\tau))^{\frac{N}{2}} d\tau \leq$$

$$(5.10) \quad \omega_N \frac{N}{2} \left(\frac{2}{\pi}\right)^{\frac{N}{2}} \int_{n^{-1}}^{t/2} \frac{\log^{\frac{N}{2}} \rho(\tau)}{\tau \rho(\tau)} d\tau,$$

where $\omega_N$ is a volume of the unit ball in $\mathbb{R}^N$. From the convergence of the integral (1.15) it follows that for some fixed value of $T$, (5.6) is satisfied for all sufficiently large $n$. This completes the proof of the $h$-irregularity of the singularity point $\mathcal{O}$.

Let us prove the $h$-regularity of the singularity point $\mathcal{O}$ by assuming that the integral (1.15) diverges. Without loss of generality, we assume that $\rho$ satisfies the

additional conditions

(5.11) $$\rho \in C^1(0, \delta),$$

(5.12) $$\frac{t\rho'(t)}{\rho(t)} \geq -L > -\frac{N}{2}$$

(5.13) $$\rho(t) > |\log t| \text{ for } 0 < t < \delta.$$

The proof of the "if" statement is based on the construction of the family of $h$-superparabolic functions $\tilde{u}_n$ with the following properties:

(1) $\tilde{u}_n(x,t) \geq 0$ in $\Omega_n$;
(2) $|1 - \tilde{u}_n(x, n^{-1})| < \frac{1}{2}$;
(3) $\forall \epsilon > 0$ there exists a number $T < \delta$ such that $\forall t_0 < T$ and for arbitrary sufficiently large $n$ we have $\tilde{u}_n(x, t_0) < \epsilon$.

Indeed, the existence of such a family implies that $u_n(x, t_0) \leq 2\tilde{u}_n(x, t_0) < 2\epsilon$, for all large $n$. passing to the limit $n \uparrow +\infty$ it follows that $u_*(x, t_0) \leq 2\epsilon$. From the maximum principle, it follows that $u_*(x,t) \leq 2\epsilon$ for all $t > t_0$. Since $\epsilon > 0$ is arbitrary, the assertion of the theorem follows.

To construct such a family $\{\tilde{u}_n\}$, we need a more precise asymptotic evaluation of $\frac{w}{h}$ in $\Omega$ as $t \downarrow 0$. We have

(5.14) $$\frac{w(\gamma, t)}{h(\gamma, t)} = -\frac{N}{2} \int_{n^{-1}}^{t} \int_{|\xi-\gamma| \leq (4\tau \log \rho(\tau))^{\frac{1}{2}}} \frac{exp\left(-\frac{|\xi-\gamma|^2}{4(t-\tau)}\right)}{(4\pi(t-\tau))^{\frac{N}{2}}} \frac{1}{\tau \rho(\tau)} \left(\frac{\tau}{t}\right)^{-\frac{N}{2}} d\xi \, d\tau.$$

We split the integral $\int_{n^{-1}}^{t} = \int_{n^{-1}}^{t\mu(t)} + \int_{t\mu(t)}^{t} =: I + J$, where $\mu(t) = k \log^{-2} \rho(t)$ and $k > 0$ is a small number at our disposal. Next, we find the asymptotics of $I$ as $t \downarrow 0$, and prove that it provides a dominating term for the asymptotic behavior of $w/h$. Since $\mu(t) \to 0$ as $t \downarrow 0$, we have $t\mu(t) \ll t/2$, and $t - \tau > t/2$ for $n^{-1} \leq \tau \leq t\mu(t)$, and $t$ sufficiently small. Therefore, we have

$$\left| -\frac{|\xi-\gamma|^2}{4(t-\tau)} \right| < \frac{|\xi-\gamma|^2}{2t} < \frac{2\tau \log \rho(\tau)}{t} < 2\mu(t) \log \rho(t\mu(t)),$$

where the last inequality follows from the monotonicity of $\tau \log \rho(\tau)$. Indeed, from the assumption (5.12) it follows that

(5.15) $$(\tau \log \rho(\tau))' = \log \rho(\tau) + \frac{\tau \rho'(\tau)}{\rho(\tau)} > 0 \text{ as } \tau \downarrow 0.$$

From (5.12) it follows that

$$\int_{t\mu(t)}^{t} \frac{\rho'(\tau)}{\rho(\tau)} d\tau = \log \rho(t) - \log \rho(t\mu(t)) \geq -\int_{t\mu(t)}^{t} \frac{L}{\tau} d\tau = L \log \mu(t),$$

which implies that

$$0 > 1 - \frac{\log \rho(t\mu(t))}{\log \rho(t)} \geq L \frac{\log \mu(t)}{\log \rho(t)} = L \frac{\log k - \log_2 \rho(t)}{\log \rho(t)} \uparrow 0,$$

and therefore,

(5.16) $$\lim_{t \downarrow 0} \frac{\log \rho(t\mu(t))}{\log \rho(t)} = 1.$$

Hence, we have
$$\left|-\frac{|\xi-\gamma|^2}{4(t-\tau)}\right| < \frac{3k}{\log\rho(t)},$$
for all sufficiently small $t$. That is to say, $\forall\,\epsilon > 0$ we can find $T > 0$ such that $\forall t < T$
$$exp\left(-\frac{|\xi-\gamma|^2}{4(t-\tau)}\right) > 1-\epsilon,$$
and therefore, we have
$$(1-\epsilon)\frac{N\omega_N}{2\pi^{\frac{N}{2}}}\int_{n^{-1}}^{t\mu(t)}\frac{\log^{\frac{N}{2}}\rho(\tau)}{\tau\rho(\tau)}\left(\frac{t}{t-\tau}\right)^{\frac{N}{2}}d\tau \leq |I|$$
$$\leq \frac{N\omega_N}{2\pi^{\frac{N}{2}}}\int_{n^{-1}}^{t\mu(t)}\frac{\log^{\frac{N}{2}}\rho(\tau)}{\tau\rho(\tau)}\left(\frac{t}{t-\tau}\right)^{\frac{N}{2}}d\tau$$
where $\omega_N$ is the volume of the unit ball. Since $n^{-1} \leq \tau \leq t\mu(t)$, we have
$$1 \leq \frac{t}{t-\tau} = 1+\frac{\tau}{t-\tau} \leq 1+\frac{\mu(t)}{1-\mu(t)} \downarrow 1, \text{ as } t\downarrow 0.$$
Therefore, $\exists\, T > 0$ such that $\forall t < T$, and for all sufficiently large $n$ we have
$$(1-\epsilon)\frac{N\omega_N}{2\pi^{\frac{N}{2}}}\int_{n^{-1}}^{t\mu(t)}\frac{\log^{\frac{N}{2}}\rho(\tau)}{\tau\rho(\tau)}d\tau \leq |I| \leq (1+\epsilon)\frac{N\omega_N}{2\pi^{\frac{N}{2}}}\int_{n^{-1}}^{t\mu(t)}\frac{\log^{\frac{N}{2}}\rho(\tau)}{\tau\rho(\tau)}d\tau$$
Hence, the following asymptotic relation is proved
$$\lim_{t\downarrow 0}\lim_{n\uparrow+\infty}\frac{I}{\frac{N\omega_N}{2\pi^{\frac{N}{2}}}\int_{n^{-1}}^{t\mu(t)}\frac{\log^{\frac{N}{2}}\rho(\tau)}{\tau\rho(\tau)}d\tau} = 1.$$
Since (1.15) is divergent, it follows that
$$(5.17)\qquad \lim_{t\downarrow 0}\lim_{n\uparrow+\infty}\frac{\frac{w(\gamma,t)}{h(\gamma,t)}}{-\frac{N\omega_N}{2\pi^{\frac{N}{2}}}\int_{n^{-1}}^{t}\frac{\log^{\frac{N}{2}}\rho(\tau)}{\tau\rho(\tau)}d\tau} = 1,$$
provided that the integrals $J$ and
$$I_1 = \int_{t\mu(t)}^{t}\frac{|\log\rho(\tau)|^{\frac{N}{2}}}{\tau\rho(\tau)}d\tau$$
remain bounded as $t\downarrow 0$. We split the integral $\int_{t\mu(t)}^{t} = \int_{t\mu(t)}^{\theta t}+\int_{\theta t}^{t} =: I_2+I_3$, with $0 < \theta < 1$ to be selected. For sufficiently small $t$ we have
$$I_3 < \int_{\theta t}^{t}\frac{d\tau}{\tau} = \log\frac{1}{\theta},\quad I_2 < \int_{t\mu(t)}^{\theta t}\frac{d\tau}{\tau\rho^{\frac{1}{2}}(\tau)} =: I_4,$$
and we still need to demonstrate the boundedness of $I_4$ as $t\downarrow 0$. This will be proved below while proving the boundedness of the integrals $I_5$ and $I_6$.

Next, we estimate $\frac{w}{h}$ inside $\Omega_n$ for small $t$. As before, we split the time integral as $\int_{n^{-1}}^{t} = \int_{\theta t}^{t} + \int_{t\mu(t)}^{\theta t} + \int_{n^{-1}}^{t\mu(t)} := I_5 + I_6 + I_7$. To estimate $I_5$, we use the identity

$$\frac{|x-\xi|^2}{4(t-\tau)} - \frac{|x-\gamma|^2}{4t} = \frac{|t(\xi-\gamma) - \tau(x-\gamma)|^2}{4t\tau(t-\tau)} - \frac{|\gamma-\xi|^2}{4\tau},$$

and derive

$$|I_5| = \frac{N}{2} \int_{\theta t}^{t} \int_{|\xi-\gamma| \leq (4\tau \log \rho(\tau))^{\frac{1}{2}}} \frac{e^{-\frac{|t(\xi-\gamma)-\tau(x-\gamma)|^2}{4t\tau(t-\tau)}}}{(4\pi(t-\tau))^{\frac{N}{2}}} \frac{1}{\tau\rho(\tau)} \left(\frac{\tau}{t}\right)^{-\frac{N}{2}} e^{\frac{|\xi-\gamma|^2}{4\tau}} d\xi\, d\tau$$

$$= \frac{N}{2} \int_{\theta t}^{t} \int_{|\xi-\gamma| \leq (4\tau \log \rho(\tau))^{\frac{1}{2}}} \frac{e^{-\frac{\left|\sqrt{\frac{t}{\tau}}(\xi-\gamma) - \sqrt{\frac{\tau}{t}}(x-\gamma)\right|^2}{4(t-\tau)}}}{(4\pi(t-\tau))^{\frac{N}{2}}} \frac{1}{\tau\rho(\tau)} \left(\frac{\tau}{t}\right)^{-\frac{N}{2}} e^{\frac{|\xi-\gamma|^2}{4\tau}} d\xi\, d\tau$$

$$\leq \frac{N}{2} \int_{\theta t}^{t} \frac{1}{\tau}\left(\frac{\tau}{t}\right)^{-\frac{N}{2}} \int_{|\xi-\gamma| \leq (4\tau \log \rho(\tau))^{\frac{1}{2}}} \frac{e^{-\frac{\left|\sqrt{\frac{t}{\tau}}(\xi-\gamma) - \sqrt{\frac{\tau}{t}}(x-\gamma)\right|^2}{4(t-\tau)}}}{(4\pi(t-\tau))^{\frac{N}{2}}} d\xi\, d\tau$$

$$= \frac{N}{2} \int_{\theta t}^{t} \frac{1}{\tau} \int_{|\xi| \leq (4t \log \rho(\tau))^{\frac{1}{2}}} \frac{e^{-\frac{\left|\xi - \sqrt{\frac{\tau}{t}}(x-\gamma)\right|^2}{4(t-\tau)}}}{(4\pi(t-\tau))^{\frac{N}{2}}} d\xi\, d\tau$$

$$\leq \frac{N}{2} \int_{\theta t}^{t} \frac{1}{\tau} d\tau \pi^{-\frac{N}{2}} \int_{\mathbb{R}^N} e^{-|z|^2} dz = \frac{N}{2} \int_{\theta t}^{t} \frac{d\tau}{\tau} = \frac{N}{2} \log \frac{1}{\theta}.$$

Let us estimate the integral

$$I_6 = -\frac{N}{2} \int_{t\mu(t)}^{\theta t} \int_{|\xi-\gamma| \leq (4\tau \log \rho(\tau))^{\frac{1}{2}}} \frac{exp\left(-\frac{|x-\xi|^2}{4(t-\tau)}\right)}{(4\pi(t-\tau))^{\frac{N}{2}}} \frac{1}{\tau\rho(\tau)} \left(\frac{\tau}{t}\right)^{-\frac{N}{2}} e^{\frac{|x-\gamma|^2}{4t}} d\xi\, d\tau$$

(5.18)
$$= -\frac{N}{2} \int_{t\mu(t)}^{\theta t} \int_{|\xi-\gamma| \leq (4\tau \log \rho(\tau))^{\frac{1}{2}}} \frac{exp\left(\frac{-\tau|x-\gamma|^2 - t|\gamma-\xi|^2 - 2t\langle x-\gamma, \gamma-\xi\rangle}{4t(t-\tau)}\right)}{(4\pi(t-\tau))^{\frac{N}{2}}} \frac{1}{\tau\rho(\tau)} \left(\frac{\tau}{t}\right)^{-\frac{N}{2}} d\xi\, d\tau.$$

Assuming $\theta < \frac{1}{2}$, from $\tau < \theta t$ it follows that $t - \tau > t/2$. Therefore, we have

$$\frac{-\tau|x-\gamma|^2 - t|\gamma-\xi|^2 - 2t\langle x-\gamma, \gamma-\xi\rangle}{4t(t-\tau)} < \frac{\langle \gamma - x, \gamma - \xi\rangle}{2(t-\tau)} < \frac{|x-\gamma||\gamma-\xi|}{t}$$

$$< \frac{(4t \log \rho(t))^{\frac{1}{2}}(4\tau \log \rho(\tau))^{\frac{1}{2}}}{t} < 4\left(\frac{\tau}{t}\right)^{\frac{1}{2}}(\log \rho(t))^{\frac{1}{2}}(\log \rho(\tau))^{\frac{1}{2}} < 4\theta^{\frac{1}{2}} \log \rho(\tau).$$

From (5.18) we deduce

$$|I_6| \leq \frac{N\omega_N}{2} \int_{t\mu(t)}^{\theta t} \frac{(4\tau \log \rho(\tau))^{\frac{N}{2}} \rho^{4\theta^{\frac{1}{2}}}(\tau)}{(2\pi t)^{\frac{N}{2}}} \frac{1}{\tau\rho(\tau)} \left(\frac{\tau}{t}\right)^{-\frac{N}{2}} dy\, d\tau$$

$$= \frac{N\omega_N}{2} \left(\frac{2}{\pi}\right)^{\frac{N}{2}} \int_{t\mu(t)}^{\theta t} \frac{\log^{\frac{N}{2}} \rho(\tau)}{\tau \rho^{1-4\theta^{\frac{1}{2}}}(\tau)} d\tau.$$

At this point, we are going to make a precise choice of the number $\theta$. Since for arbitrary $\gamma > 0$

$$\frac{\log^{\frac{N}{2}} \rho(\tau)}{\rho^\gamma(t)} \to 0 \text{ as } t \downarrow 0,$$

we can reduce the boundedness of $|I_6|$ to the boundedness of $|I_4|$ if we choose $\theta$ such that

$$1 - 4\theta^{\frac{1}{2}} > \frac{1}{2} \text{ or } \theta < \frac{1}{64}.$$

We fix the value $\theta = \frac{1}{65}$. Therefore, for sufficiently small $t$ we have

$$|I_6| \leq \frac{N\omega_N}{2} \left(\frac{2}{\pi}\right)^{\frac{N}{2}} \int_{t\mu(t)}^{\frac{1}{65}t} \frac{d\tau}{\tau\rho^{\frac{1}{2}}(\tau)} d\tau.$$

By using assumption (5.13) we have

$$\int_{t\mu(t)}^{\frac{t}{65}} \frac{d\tau}{\tau\rho^{\frac{1}{2}}(\tau)} d\tau \leq C \int_{t\mu(t)}^{\frac{t}{65}} \frac{d\tau}{\tau |\log \tau|^{\frac{1}{2}}} d\tau = 2\left[|\log t\mu(t)|^{\frac{1}{2}} - \left|\log \frac{t}{65}\right|^{\frac{1}{2}}\right]$$

$$= 2\frac{|\log t\mu(t)| - \left|\log \frac{t}{65}\right|}{\left|\log \frac{t}{65}\right|^{\frac{1}{2}} + |\log t\mu(t)|^{\frac{1}{2}}} = 2\frac{|\log \mu(t)| - \log 65}{\left|\log \frac{t}{65}\right|^{\frac{1}{2}} + |\log t\mu(t)|^{\frac{1}{2}}}$$

(5.19) $$2\frac{|\log \mu(t)| - \log 65}{\left|\log \frac{t}{65}\right|^{\frac{1}{2}}\left[1 + \left|\frac{\log t\mu(t)}{\log t/65}\right|^{\frac{1}{2}}\right]} \leq \frac{|\log \mu(t)|}{\left|\log \frac{t}{65}\right|^{\frac{1}{2}}}.$$

By using a l'Hopital's rule and (5.12) we have

(5.20) $$\lim_{t\to 0} \frac{\log^2 \mu(t)}{-\log t/65} = -4 \lim_{t\to 0} \frac{2\log\log \rho(t) - \log k}{\log \rho(t)} \frac{t\rho'(t)}{\rho(t)} = 0,$$

and therefore, from (5.19) we deduce that

(5.21) $$\lim_{t\to 0} \int_{t\mu(t)}^{\frac{t}{65}} \frac{d\tau}{\tau\rho^{\frac{1}{2}}(\tau)} d\tau = 0.$$

Finally, we estimate the integral

$$I_7 = -\frac{N}{2} \int_{n^{-1}}^{t\mu(t)} \int_{|\xi-\gamma|\leq(4\tau\log\rho(\tau))^{\frac{1}{2}}} \frac{exp\left(-\frac{|x-\xi|^2}{4(t-\tau)}\right)}{(4\pi(t-\tau))^{\frac{N}{2}}} \frac{1}{\tau\rho(\tau)} \left(\frac{\tau}{t}\right)^{-\frac{N}{2}} e^{\frac{|x-\gamma|^2}{4t}} d\xi\, d\tau$$

$$= -\frac{N}{2} \int_{n^{-1}}^{t\mu(t)} \int_{|\xi-\gamma|\leq(4\tau\log\rho(\tau))^{\frac{1}{2}}} \frac{exp\left(\frac{-\tau|x-\gamma|^2-t|\gamma-\xi|^2-2t\langle x-\gamma,\gamma-\xi\rangle}{4t(t-\tau)}\right)}{(4\pi(t-\tau))^{\frac{N}{2}}} \frac{1}{\tau\rho(\tau)} \left(\frac{\tau}{t}\right)^{-\frac{N}{2}} d\xi\, d\tau.$$

Our goal is to demonstrate that its asymptotics as $t \downarrow 0$ coincides with the asymptotics of the corresponding integral $I$ in the expression of $\frac{w(\gamma,t)}{h(\gamma,t)}$. To prove that we need to demonstrate that the term $exp\left(\frac{-\tau|x-\gamma|^2-2t\langle x-\gamma,\gamma-\xi\rangle}{4t(t-\tau)}\right)$ is close to 1 for small $t$. We have

$$\left|\frac{-\tau|x-\gamma|^2 - 2t\langle x-\gamma,\gamma-\xi\rangle}{4t(t-\tau)}\right| \leq \frac{\tau|x-\gamma|^2}{4t(t-\tau)} + \frac{|x-\gamma||\gamma-\xi|}{2(t-\tau)}$$

Since $t-\tau > \frac{t}{2}$, we have

$$\frac{\tau|x-\gamma|^2}{4t(t-\tau)} \leq \frac{2\tau\log\rho(t)}{t} \leq \frac{2t\mu(t)\log\rho(t)}{t} = \frac{2k}{\log\rho(t)} \to 0, \text{ as } t \downarrow 0.$$

Using (5.15), (5.16), we also deduce that for all sufficiently small $t$

$$\frac{|x-\gamma||\gamma-\xi|}{2(t-\tau)} \leq \frac{|x-\gamma||\gamma-\xi|}{t} \leq 4\left(\frac{\tau\log\rho(\tau)\log\rho(t)}{t}\right)^{\frac{1}{2}}$$

$$\leq 4(\mu(t)\log\rho(t\mu(t))\log\rho(t))^{\frac{1}{2}} = 4k^{\frac{1}{2}}\left(\frac{\log\rho(t\mu(t))}{\log\rho(t)}\right)^{\frac{1}{2}} < 5k^{\frac{1}{2}},$$

and the right-hand side will be sufficiently small if the parameter $k$ at our disposal is chosen small enough.

Hence, we proved that for arbitrary $\epsilon > 0$ there exists $T > 0$ such that for any fixed $t < T$ and for all large $n$

$$\left|\frac{\frac{w(x,t)}{h(x,t)}}{\frac{w(\gamma,t)}{h(\gamma,t)}} - 1\right| < \epsilon \text{ for } (x,t) \in \Omega_n.$$

Otherwise speaking,

$$\lim_{t\downarrow 0} \lim_{n\uparrow+\infty} \frac{\frac{w(x,t)}{h(x,t)}}{\frac{w(\gamma,t)}{h(\gamma,t)}} = 1, \text{ uniformly for all } x \text{ with } (x,t) \in \Omega_n.$$

Consider a function

$$\tilde{u}_n(x,t) = \frac{\frac{v(x,t)+w_1(x,t)}{h(x,t)}}{\sup_{\Omega_n}\left|\frac{w_1(x,t)}{h(x,t)}\right|} + 1,$$

where

$$w_1(x,t) = \frac{N-2L}{N}w(x,t).$$

From (5.4), (5.5) and (5.12) it follows that $\mathcal{H}(\tilde{u}_n h) \geq 0$ for small $t$. Therefore, $\tilde{u}_n$ is $h$-superparabolic function for sufficiently small t. We can check that it satisfies the required conditions (1)-(3).

(1) We have $|\frac{v}{h}| < 1$ in $\Omega_n$, $w_1(x, n^{-1}) = 0$, and $\frac{w_1(y,t)}{h(y,t)} \to -\infty$ as $n \uparrow +\infty$, $t \downarrow 0$. Therefore,

(5.22) $$\tilde{u}_n(x, n^{-1}) \to 1 \text{ as } n \uparrow +\infty \text{ uniformly for } x.$$

(2) We have

(5.23) $$\lim_{t \downarrow 0} \lim_{n \uparrow +\infty} \frac{\frac{w_1(\gamma,t)}{h(\gamma,t)}}{-\frac{(N-2L)\omega_N}{2\pi^{\frac{N}{2}}} \int_{n^{-1}}^{t} \frac{\log^{\frac{N}{2}} \rho(\tau)}{\tau \rho(\tau)} d\tau} = 1,$$

(5.24) $$\lim_{n \uparrow +\infty} \int_{n^{-1}}^{t} \frac{\log^{\frac{N}{2}} \rho(\tau)}{\tau \rho(\tau)} d\tau = +\infty,$$

(5.25) $$\lim_{t \downarrow 0} \lim_{n \uparrow +\infty} \frac{\frac{w_1(x,t)}{h(x,t)}}{\frac{w_1(\gamma,t)}{h(\gamma,t)}} = 1, \text{ uniformly for all } x \text{ with } (x,t) \in \Omega_n.$$

From these conditions it follows that $\forall \epsilon > 0$ there exists a number $T < \delta$ such that $\forall 0 < t_0 < T$ and for arbitrary sufficiently large $n$ we have $\tilde{u}_n(x, t_0) < \epsilon$.

(3) $\tilde{u}_n(x,t) \geq 0$ in $\Omega_n$, since $\frac{v}{h} \geq 0$ and $\frac{w_1}{h} \leq 0$ in $\Omega_n$.

Hence we proved that the divergence of the integral (1.15) implies the $h$-regularity of $\mathcal{O}$ for $\Omega$ provided that additional assumptions (5.11)-(5.13) are satisfied. Note that the assumptions (5.11)-(5.13) are satisfied for all functions in (1.16). Therefore, we completed the proof of $h$-regularity of $\mathcal{O}$ and removability of the fundamental singularity for domains (1.18) with $\epsilon \leq 0$.

To complete the proof we only need to demonstrate that the assumptions (5.11)-(5.13) can be removed. Differentiability assumption (5.11) can be removed as before with the only difference that we select $\rho_1 \in C^1(0, \delta)$ which satisfy

(5.26) $$\rho(t) < \rho_1(t) < 2\rho(t), \ 0 < t < \delta,$$

and consider $\Omega_1$ as in (5.2). From (5.37) it follows that the integral (1.15) is divergent for $\rho_1$, and $\Omega_1$ contains $\Omega \cap \{t < \delta\}$. Therefore, $h$-regularity of $\mathcal{O}$ for $\Omega_1$ implies the $h$-regularity of $\mathcal{O}$ for $\Omega$.

To remove assumption (5.13), assume on the contrary that there are arbitrarily small values of $t$ such that $\rho(t) \leq |\log t|$. Consider a function

$$\rho_1(t) = \max(\rho(t); |\log t|)$$

Clearly, the $h$-regularity of $\mathcal{O}$ for $\Omega_1$ implies the $h$-regularity of $\mathcal{O}$ for $\Omega$. Hence, we only need to demonstrate that the integral (1.15) is divergent for $\rho_1$. In view of our assumption we can choose the sequence $\{t_k\}$ with the following properties:

(5.27) $$\begin{cases} t_k > t_{k+1} > \cdots > 0, \ t_k \downarrow 0 \\ \rho(t_k) = |\log t_k| \\ \frac{\log t_{k+1}}{\log t_k} \geq 2 \end{cases}$$

Let us define a function

$$\tilde{\rho}(t) = |\log t_{k+1}|, \text{ for } t_{k+1} \leq t < t_k.$$

We have

(5.28) $$\tilde{\rho}(t) \geq \rho_1(t), \ 0 < t < t_1.$$

By the last relation of (5.27) we have

$$\int\limits_{t_{k+1}}^{t_k} \frac{dt}{t\tilde{\rho}(t)} = 1 - \frac{\log t_k}{\log t_{k+1}} \geq \frac{1}{2},$$

which implies that

(5.29) $$\int_{0+} \frac{dt}{t\tilde{\rho}(t)} = +\infty.$$

From (5.28),(5.29) it follows that

(5.30) $$\int_{0+} \frac{dt}{t\rho_1(t)} = +\infty,$$

and therefore, the integral (1.15) is also divergent for $\rho_1$. Hence, assumption (5.13) is removable.

Next, we are going to demonstrate that the assumption (5.12) is also removable. Assume that the given function $\rho \in C^1(0, \delta)$ has a divergent integral (1.15), but doesn't satisfy the condition (5.12). Consider a one-parameter family of curves

(5.31) $$\rho_C(t) = |\log Ct|^3, \ C > 0, \ 0 < t < C^{-1}.$$

Note that for any curve $\rho_C$ the integral (1.15) is convergent, and since $t^{-1} \log^{\frac{N}{2}} t$ is monotonically decreasing function for large $t$, the convergence of (1.15) holds for any function satisfying $\rho \geq \rho_C$ near 0. Therefore, there are arbitrarily small values of $t$ such that $\rho(t) < \rho_C(t)$. Since $\rho_C(C^{-1}) = 0$, a graph of the given function $\rho$ with divergent integral (1.15) must intersect all curves $\rho_C(t)$ with $C \geq \delta^{-1}$. For any point $(\rho(t), t)$ on the positive quarter plane there exists a unique value

(5.32) $$C = C(t) = t^{-1} e^{-\rho^{\frac{1}{3}}(t)}$$

such that $\rho_C(t)$ passes through the point $(\rho(t), t)$. Clearly, $tC(t) \to 0$ as $t \to 0$; and there exists a sequence $t_n \downarrow 0$ such that $C(t_n) \uparrow +\infty$ as $n \to +\infty$. We define a set

(5.33) $$M = \{t \in (0, \delta] : C(t) = C_1(t)\},$$

where $C_1(t) = \max\limits_{t \leq \tau \leq \delta} C(\tau)$. Denote by $\overline{M}^c$ the complement of $\overline{M}$. Since $\overline{M}^c$ is an open set, we have

(5.34) $$\overline{M}^c = \cup_n (t_{2n}, t_{2n-1}).$$

From the definition of the set $M$ it follows that $C(t)$ satisfies the following properties:

(5.35) $$\begin{cases} C(t) \text{ is a decreasing function for } t \in \overline{M} \\ C(t_{2n}) = C(t_{2n-1}) \text{ for } (t_{2n}, t_{2n-1}) \subset \overline{M}^c \\ C(t) \uparrow +\infty \text{ as } \overline{M} \ni t \downarrow 0. \end{cases}$$

Indeed, if we take $t', t'' \in M$ with $t' > t''$, then we have

$$C(t'') = C_1(t'') \geq C_1(t') = C(t').$$

On the other hand, if $t', t'' \in \overline{M}$ the same conclusion follows from the continuity of the function $C(t)$.

To prove the second assertion, first note that since $t_{2n}, t_{2n-1} \in \overline{M}$, we have $C_1(t_{2n}) = C(t_{2n})$ and $C_1(t_{2n-1}) = C(t_{2n-1})$. Therefore, assuming that $C(t_{2n}) \neq C(t_{2n-1})$ would imply that $C_1(t_{2n}) > C_1(t_{2n-1})$. Since $C_1$ is a continuous function for some $\epsilon \in (0, t_{2n-1}-t_{2n})$ we have $C_1(t_{2n-1}) < C_1(t_{2n}+\epsilon)$. Let $C_1(t_{2n}+\epsilon) = C(\theta)$. Obviously, we must have $\theta \in [t_{2n}+\epsilon, t_{2n-1})$ and $C_1(\theta) = C(\theta)$. However, this is a contradiction with the fact that $(t_{2n}, t_{2n-1}) \subset \overline{M}^c$, which proves the second assertion of (5.35).

To prove the third assertion of (5.35), recall that there is a sequence $t_n \downarrow 0$ with $C(t_n) \uparrow +\infty$. Assume that for some $n$ we have $t_n \notin \overline{M}$. Then there exists $k_n$ such that $t_n \in (t_{k_n}, t_{k_n+1}) \subset \overline{M}^c$. According to the second assertion of (5.35) we have

$$C(t_n) \leq C_1(t_n) \leq C_1(t_{k_n}) = C(t_{k_n}) = C(t_{k_n+1})$$

Consider a sequence

$$(5.36) \qquad \tilde{t}_n = \begin{cases} t_n, & \text{if } t_n \in \overline{M}, \\ t_{k_n}, & \text{if } t_n \in (t_{k_n}, t_{k_n+1}) \subset \overline{M}^c. \end{cases}$$

Clearly, we have $\tilde{t}_n \in \overline{M}$; $\tilde{t}_n \downarrow 0$ and $C(\tilde{t}_n) \uparrow +\infty$ as $n \uparrow +\infty$. From the monotonicity of $C(t)$ on $\overline{M}$ we easily deduce the third assertion of (5.35).

Let us define now a new function $\rho_1$ with the following properties:

(1) $\rho_1(t) = \rho(t)$ for $t \in \overline{M}$,

(2) $\rho_1(t) = |\log(C(t_{2n})t)|^3$ for $t \in (t_{2n}, t_{2n-1}) \subset \overline{M}^c$.

Note that equivalently we can define $\rho_1$ as

$$(5.37) \qquad \rho_1(t) = |\log(C_1(t)t)|^3, \ 0 < t \leq \delta.$$

In fact, function $C(t)$ defined for the function $\rho_1$ via (5.32) coincides with $C_1(t)$, i.e.

$$C(t) = \max_{t \leq \tau \leq \delta} C(\tau).$$

The new function is continuous and satisfy $\rho_1(t) \geq \rho(t)$ everywhere, with $\rho_1(t) \neq \rho(t)$ on intervals $(t_{2n}, t_{2n-1})$.

It can be easily seen that the function $\rho$ satisfies the desired property (5.12) on $\overline{M}$. Indeed, for a function $\rho_C$ we have

$$\frac{t\rho'_C(t)}{\rho_C(t)} = \frac{-3}{|\log Ct|} \to 0, \text{ as } Ct \to 0.$$

Since $C(t)$ is monotonically decreasing on $\overline{M}$, we have

$$|\rho'(t)| \leq |\rho'_C(t)| = \left| \frac{3\log^2(Ct)}{t} \right|, \ t \in \overline{M}$$

provided that $C = C(t)$ is chosen as in (5.32). Therefore, we have

$$(5.38) \qquad \left| \frac{t\rho'(t)}{\rho(t)} \right| = \left| \frac{3}{\log(C(t)t)} \right|, \ t \in \overline{M}.$$

Since $C(t)t \to 0$ as $t \downarrow 0$, the right-hand side is arbitrarily small for small $t$, and clearly $\rho$ satisfies (5.12) on $\overline{M}$.

Our next goal is to demonstrate that

$$(5.39) \qquad \int_{\overline{M}^c} \frac{\log^{\frac{N}{2}} \rho(t)}{t\rho(t)} dt < +\infty$$

Since $\rho \geq \rho_1$ it is satisfactory to demonstrate that the integral (5.39) is convergent for $\rho_1$. By using the property (2) of $\rho_1$ we have

$$(5.40) \qquad \int_{\overline{M}^c} \frac{\log^{\frac{N}{2}} \rho_1(t)}{t\rho_1(t)} \, dt = \sum_n \int_{t_{2n}}^{t_{2n-1}} \frac{3^{\frac{N}{2}} \log^{\frac{N}{2}} |\log C_n t|}{t |\log C_n t|^3} \, dt,$$

where we use a notation $C_n := C(t_{2n}) \equiv C(t_{2n-1})$. Since $C_n t$ is arbitrarily small for $t$ sufficiently small, we can use the following inequality for all large n:

$$(5.41) \qquad \frac{3^{\frac{N}{2}} \log^{\frac{N}{2}} |\log C_n t|}{|\log C_n t|} \leq 1$$

Using (5.41), from (5.40) we deduce
(5.42)
$$\int_{\overline{M}^c} \frac{\log^{\frac{N}{2}} \rho_1(t)}{t\rho_1(t)} \, dt \leq \sum_n \int_{t_{2n}}^{t_{2n-1}} \frac{dt}{t |\log C_n t|^2} \, dt = \sum_n \frac{1}{|\log C_n t_{2n-1}|} - \frac{1}{|\log C_n t_{2n}|}$$

Since
$$\rho(t_{2n-1}) = |\log C_n t_{2n-1}|^3, \rho(t_{2n}) = |\log C_n t_{2n}|^3,$$
and both $\rho(t)$ and $|\log t|^3$ are decreasing functions, we deduce that

$$(5.43) \qquad C_n t_{2n-1} \geq C_n t_{2n} \geq C_{n+1} t_{2n+1}, \ n = 1, 2, ...$$

Therefore we also have

$$(5.44) \qquad \frac{1}{|\log C_n t_{2n-1}|} \geq \frac{1}{|\log C_n t_{2n}|} \geq \frac{1}{|\log C_{n+1} t_{2n+1}|}, \ n = 1, 2, ...$$

Hence, the series (5.42) is a telescoping series and therefore integrals (5.40) and (5.39) are convergent integrals. Since $\rho_1 \leq \rho$ on $\overline{M}^c$, it follows that the integral (5.39) is convergent for $\rho$. On the other side, since the integral (1.15) is divergent it follows that

$$(5.45) \qquad \int_{\overline{M}} \frac{\log^{\frac{N}{2}} \rho(t)}{t\rho(t)} \, dt = +\infty.$$

Now we pursue the identical proof given above under the assumptions (5.11)-(5.13) by replacing the function (5.8) with the following one:

$$(5.46) \qquad \tilde{w}(x, t) = -\frac{N}{2} \int_{n^{-1}}^{t} \int_{|\xi-\gamma| \leq (4\tau \log \rho(\tau))^{\frac{1}{2}}} \frac{exp\left(-\frac{|x-\xi|^2}{4(t-\tau)}\right)}{(4\pi(t-\tau))^{\frac{N}{2}}} \frac{(4\pi\tau)^{-\frac{N}{2}}}{\tau \tilde{\rho}(\tau)} \, d\xi \, d\tau,$$

where the function $\tilde{\rho}$ in the integrand is chosen with the following properties:

(1) $\tilde{\rho}(t) = \rho(t)$ for $t \in \overline{M}$;

(2) In the complementary set $\overline{M}^c$ the function $\tilde{\rho}$ is chosen as sufficiently large continuous function with the property that all the estimations of the function $w/h$ from (5.14) in the proof given above remain valid when the function $\rho$ in the integrand is replaced with $\tilde{\rho}$.

(3) $\tilde{\rho}(t) \geq \rho(t), \ 0 < t \leq \delta.$

To estimate the function $\tilde{w}/h$ in $\Omega$ as $t \downarrow 0$, the time integral in $(n^{-1}, t)$ is split into two parts over $\overline{M} \cap (n^{-1}, t)$ and $\overline{M}^c \cap (n^{-1}, t)$. The estimation of the first one is identical to the presented proof, for the assumptions (5.11)-(5.13) are satisfied on $\overline{M}$. Due to property (2) of the function $\tilde{\rho}$ the second integral remains bounded and accordingly does not affect the leading asymptotic of $\tilde{w}/h$ given via the divergent integral (5.45). Precisely, we establish (5.17),(5.22)-(5.25), where the integral term $\int_{n^{-1}}^{t}$ in expressions (5.23),(5.24) is replaced with $\int_{\overline{M} \cap (n^{-1}, t)}$. This completes the proof of the $h$-regularity of $\mathcal{O}$ without assumptions (5.11)-(5.13). □

Theorem 1.3 follows from the Lemma 4.1 (iv) and the mapping (4.1).

5.2. **Proof of Theorems 3.1 and 3.2.** According to the Definitions 2.4, 2.5, and the formulae (2.11), (2.16), Theorem 1.2 and 1.3 are equivalent to the claim $(1) \Leftrightarrow (4)$ of Theorems 3.1 and 3.2 respectively..

The equivalence $(1) \Leftrightarrow (2)$ follows from the formulae (2.6) and (2.8).

The equivalence $(2) \Leftrightarrow (3)$ in Theorem 3.2 with $\gamma = 0$ is proved in [1] (see Lemma 2.3, p. 472). Applying Lemma 4.1 (iv), the equivalence $(2) \Leftrightarrow (3)$ in Theorem 3.1 with $\gamma = 0$ follows. Applying the translation $x \mapsto x + \gamma$, the equivalence of $(2) \Leftrightarrow (3)$ in Theorem 3.1 with $\gamma \neq 0$ easily follows. Applying Lemma 4.1 (iv) again, the equivalence $(2) \Leftrightarrow (3)$ in Theorem 3.2 with $\gamma \neq 0$ follows. □

## References

1. U.G. Abdulla, Wiener's Criterion at $\infty$ for the Heat Equation, *Advances in Differential Equations*, **13**, 5-6(2008), 457-488.
2. U.G. Abdulla, First boundary value problem for the difusion equation. I. Iterated logarithm test for the boundary regularity and solvability, *SIAM J. Math. Anal.* **34**(2003), 1422-1434.
3. U.G. Abdulla, Multidimensional Kolmogorov-Petrovsky test for the boundary regularity and irregularity of solutions to the heat equation, *Boundary Value Problems*, **2005** (2005), 181-199.
4. U. G. Abdulla, Kolmogorov Problem for the Heat Equation and its Probabilistic Counterpart, *Nonlinear Analysis*, **63**, 2005, 712-724.
5. U.G. Abdulla, Wiener's Criterion for the Unique Solvability of the Dirichlet Problem in Arbitrary Open Sets with Non-Compact Boundaries, *Nonlinear Analysis,* **67**,2 (2007), 563-578.
6. U.G. Abdulla, Regularity of $\infty$ for elliptic equations with measurable coefficients and its consequences, *Discrete and Continuous Dynamical Systems*, **32**,10 (2012), 3379-3397.
7. U.G. Abdulla, Removability of the Logarithmic Singularity for the Elliptic PDEs with Measurable Coefficients and its Consequences, *Calculus of Variations and Partial Differential Equations*, **57:157**, 2018.
8. D.G. Aronson, Removable singularities for linear parabolic equations, *Archive for Rational Mechanics and Analysis,* **17**, 1(1964), 79-84.
9. D.G. Aronson, Isolated singularities of solutions of second order parabolic equations, *Archive for Rational Mechanics and Analysis,* **19**, 3(1965), 231-238.
10. D.G. Aronson, Non-negative solutions of linear parabolic equations, *Annali della Scuola Normale Superiore di Pisa,* **22**, 4(1968), 607-694.
11. H. Bauer, Harmonische Räume und ihre Potentialtheorie, Lecture Notes in Mathematics, **22**, Springer, 1966.
12. M. Brelot, Sur le role du point a l'infini dans la theorie des fonctions harmoniques, *Ann. Sc. Ecolo Norm. Sup.,* **61** (1944), 301-332.
13. M. Brelot, *On Topologies and Boundaries in Potential Theory*, Lecture Notes in Mathematics, **175**, Springer-Verlag, 1971. 5555
14. G. Choquet, Theory of capacities, *Ann. Inst. Fourier Grenoble,* **5**(1955), 131-295.
15. J.L. Doob, *Classical Potential Theory and its Probabilistic Counterpart*, Springer, 1984.
16. L.C. Evans and R. Gariepy, Wiener's criterion for the heat equation, *Archive for Rational Mechanics and Analysis,* **78** (1982), 293-314.

17. E.Fabes, D.Jerison and C.Kenig, The Wiener test for degenerate elliptic equations, *Annales de l'Institut Fourier*, **32**,3(1982),151-182.
18. E. Fabes, N. Garofalo and E. Lanconelli, Wiener's criterion for divergence form parabolic operators with $C^1$-Dini continuous coefficients, *Duke Math J.*, **59**, 1(1989), 191-232.
19. R.Gariepy and W.P.Ziemer, A regularity condition at the boundary for solutions of quasilinear elliptic equations, *Archive for Rational Mechanics and Analysis*, **67**(1977), 25-39.
20. N. Garofalo and E. Lanconelli, Wiener's criterion for parabolic equations with variable coefficients and its consequences, *Trans. Amer. Math. Soc.*, **308**, 2(1988), 811-836.
21. J Heinonen, T Kilpeläinen, O Martio, *Nonlinear Potential Theory of Degenerate Elliptic Equations*, Clarendon Press, 1993.
22. K. Ito and H.P. McKean,Jr., *Diffusion Processes and Their Sample Paths*. Springer, 1996.
23. T.Kilpeläinen and J. Maly, The Wiener test and potential estimates for quasilinear elliptic equations, *Acta Mathematica*, **172**(1994),137-161.
24. E. Lanconelli, Sul problema di Dirichlet per l'equazione del calore, *Ann. Math. Pura. Appl.*, **97** (1973), 83-114.
25. W. Littman, G. Stampacchia and H.F. Weinberger, Regular Points for Elliptic Equations with Discontinuous Coefficients, *Ann. Sc. Norm. Super. Pisa Cl. Sci.*, **17**(3), (1963), 43–77.
26. V.G.Maz'ya, On the continuity at a boundary point of solutions of quasi-linear elliptic equations, *Vestnik Leningrad University: Mathematics*, **3**(1976),225-242.
27. I.G. Petrovsky, Zur ersten Randwertaufgabe der Wärmeleitungsgleichung, Compositio Math. 1(1935), 383-419.
28. J.Serrin and H.F.Weinberger, Isolated singularities of solutions of linear elliptic equations, *Amer. J. Math.*, **72** (1966), 258-272.
29. A.N. Tikhonov, Uniqueness theorems for the heat equation, *Mathematichesky Sbornik, 42(2)*, 1935, 199-216.
30. N.A. Watson, Introduction to Heat Potential Theory, American Mathematical Society, 2012.
31. N. Wiener, Certain Notions in Potential Theory, *J. Math. Phys.*, **3**, (1924), 24–51.
32. N. Wiener, The Dirichlet Problem, *J. Math. Phys.*, **3**, (1924), 127–146.

ANALYSIS & PDE UNIT, OKINAWA INSTITUTE OF SCIENCE AND TECHNOLOGY, OKINAWA, JAPAN